%% file: main.tex
\documentclass[a4paper,review]{cas-sc}
\usepackage[section]{placeins} % For '\input{file}'
\input{preamble}

\begin{document}

\title[mode=title]{Positivity-Preserving Entropy-Based Adaptive Filtering for Discontinuous Spectral Element Methods}
\shorttitle{}
\shortauthors{T. Dzanic \textit{et al.}}

\author[1]{T. Dzanic}[orcid=0000-0003-3791-1134]
\cormark[1]
\cortext[cor1]{Corresponding author}
\ead{tdzanic@tamu.edu}
\author[1]{F. D. Witherden}[orcid=0000-0003-2343-412X]

\address[1]{Department of Ocean Engineering, Texas A\&M University, College Station, TX 77843}

\begin{abstract}
In this work, we present a positivity-preserving entropy-based adaptive filtering method for shock capturing in discontinuous spectral element methods. By adapting the filter strength to enforce positivity and a local discrete minimum entropy principle, the resulting approach can robustly resolve strong discontinuities with sub-element resolution, does not require problem-dependent parameter tuning, and can be easily implemented on general unstructured meshes with relatively low computational cost. The efficacy of the approach is shown in numerical experiments on hyperbolic and mixed hyperbolic-parabolic conservation laws such as the Euler and Navier-Stokes equations for
problems including extreme shocks, shock-vortex interactions, and complex compressible
turbulent flows.
\end{abstract}


%\linenumbers

% ================================================================================

\maketitle

\input{introduction}

\input{preliminaries}
\input{methodology}

\input{implementation}
\input{results}
\input{conclusions}

\section*{Acknowledgements}
\label{sec:ack}
This work was supported in part by the U.S. Air Force Office of Scientific Research via grant FA9550-21-1-0190 ("Enabling next-generation heterogeneous computing for massively parallel high-order compressible CFD") of the Defense University Research Instrumentation Program (DURIP) under the direction of Dr. Fariba Fahroo.

\bibliographystyle{unsrtnat}
\bibliography{reference}

\appendix
\input{app1}


% Show the list of todo's in the document.  Needed to avoid stupid warnings/errors when using the todo package
%\todos

\end{document}

%% file: preamble.tex
\usepackage{amssymb}
\usepackage{amsthm}
\usepackage{amsmath}
\usepackage{upgreek}
\usepackage{pdflscape}
\usepackage{listings}
\usepackage{multirow}
\usepackage[percent]{overpic}
\usepackage{multicol}
\usepackage{color}
\usepackage{stmaryrd}
\usepackage{algorithm}
\usepackage{algorithmic}
\usepackage{xfrac}
\usepackage[capitalise]{cleveref}
\usepackage{booktabs}
\usepackage{soul}
\usepackage[section]{placeins}
\usepackage{calc}
\usepackage{xcolor}
\usepackage[titletoc]{appendix}
\usepackage{siunitx}
\usepackage[sort&compress,square,numbers]{natbib}

\usepackage{graphicx}
\usepackage{graphics}
\usepackage{wrapfig}
\usepackage{float}
\usepackage{subfig}
\usepackage[percent]{overpic}
\usepackage{varwidth}
\usepackage{tikz}
\usepackage{pgfplots}
\usepackage{adjustbox}
\usetikzlibrary{arrows,matrix,positioning,fit}
\usetikzlibrary{shapes,positioning}
\usetikzlibrary{backgrounds}
\usepackage{tikz-layers}
\usepgfplotslibrary{fillbetween}
\usetikzlibrary{intersections}
\pgfplotsset{compat=1.14}
\usepgfplotslibrary{colorbrewer}
\usepgfplotslibrary{patchplots}
\usepgfplotslibrary[colorbrewer]
\usetikzlibrary{pgfplots.colorbrewer}
\usetikzlibrary[pgfplots.colorbrewer]
\usepgfplotslibrary{units}
\usetikzlibrary{spy}
\usepackage{pgfplotstable}
\usepackage{arrayjobx}
\graphicspath{ {./figs/} }
\usetikzlibrary{external}

\newlength\myheight
\newlength\mydepth
\settototalheight\myheight{Xygp}
\newcommand*\inlinegraphics[1]{%
  \settototalheight\myheight{Xygp}%
  \settodepth\mydepth{Xygp}%
  \raisebox{-\mydepth}{\includegraphics[height=\myheight]{#1}}%
}
\newcommand\orcid[1]{\href{https://orcid.org/#1}{\inlinegraphics{orcid_16x16.png}}}

\makeatletter
\def\BState{\State\hskip-\ALG@thistlm}
\makeatother

\newdefinition{definition}{Definition}[section]
\newtheorem*{remark}{Remark}

\newcommand\pxvar[2]{\partial_{#2} #1}

\newcommand{\shalf}{{\sfrac{1}{2}}}

%% file: introduction.tex
 %!TEX root = ./main.tex
\section{Introduction}
\label{sec:intro}

The use of discontinuous spectral element methods (DSEM) has grown in prevalence over the years due to their ability to achieve high-order accuracy while retaining geometric flexibility and a compact data structure suited for massively-parallel computing. As a result, DSEM offer many advantages for simulations of complex problems. However, their robustness is severely reduced for problems that exhibit discontinuities as the presence of spurious oscillations due to Gibbs phenomena can result in nonphysical solutions or the failure of the scheme altogether. Consequently, this lack of reliability is one of the limitations preventing the widespread adoption of these methods in the industry. To extend the use of DSEM to a wider variety of problems, various stabilization techniques have been proposed to increase the robustness of these schemes in the vicinity of discontinuities. A common goal of these shock capturing methods is to suppress numerical instabilities in the vicinity of a shock without degrading the accuracy of the underlying numerical scheme in regions where the solution is smooth.

These stabilization techniques can be broadly categorized as artificial viscosity, limiting, stencil modification, or filtering methods. The most ubiquitous approach is the addition of artificial viscosity to explicitly introduce numerical dissipation in the vicinity of a shock, a method pioneered by \citet{VonNeumann1950} with various advancements and alternative approaches over the decades \citep{Guermond2011, Barter2007, Persson2006, Nazarov2012, Dzanic2021}.
Limiting methods have also shown promise \citep{Boris1997, Guermond2019}, with the premise of these schemes generally relying on a combination of a constraint-satisfying low-order scheme and a constraint-violating high-order scheme. Additionally, modifying the numerical stencil to alleviate the issues of Gibbs phenomena has shown to be accurate and robust in other numerical settings \citep{Jiang1999, Shu1988, Shu1998, Harten1987}, but the use of these methods is not very prevalent in the context of DSEM, primarily due to the incompatibility of the approach with the compact data structure, local compute, and geometric flexibility afforded by DSEM. In contrast, filtering \citep{Glaubitz2017, Hesthaven2008, Panourgias2016}, where spurious oscillations are removed by reducing high-frequency modes in the solution, is particularly attractive as a shock capturing method for DSEM as the filter can sharply resolve discontinuities with minimal computational cost and without sacrificing the computational efficiency of the underlying numerical scheme. 

However, a typical drawback of many of these approaches is that they (1) do not necessarily \textit{guarantee} that physical constraints on the solution are satisfied, which may lead to the failure of the scheme, (2) have free parameters which can require problem- and mesh-dependent tuning, a cost that cannot be afforded for large scale-resolving simulations, and/or (3) are not easily and efficiently implemented in the context of explicit DSEM on modern computing architectures. These issues motivate the development of shock capturing approaches for DSEM that can guarantee certain physical constraints are satisfied without requiring problem-dependent tunable parameters. In this work, we present an adaptive filtering approach for shock capturing in nodal DSEM to address these issues. By formulating physical constraints such as positivity and a local minimum entropy principle as constraints on the discrete solution, the filter strength is computed via a simple scalar optimization problem requiring only element-local information. Under some basic assumptions on the properties of the numerical scheme, the filtered solution is guaranteed to satisfy these constraints, resulting in an efficient and robust method for resolving discontinuous features without the use of problem-dependent tunable parameters. Furthermore, the proposed filtering approach does not appreciably degrade the efficiency and accuracy of the standard DSEM approach for smooth solutions.

The remainder of this paper is organized as follows. Preliminaries regarding the methodology are presented in \cref{sec:prelim}. The proposed filtering approach is then presented in \cref{sec:methodology} followed by numerical implementation details in \cref{sec:implementation}. The results of numerical experiments on the Euler and Navier--Stokes equations are then given in \cref{sec:results}, followed by conclusions drawn in \cref{sec:conclusion}.

%% file: preliminaries.tex
\section{Preliminaries}\label{sec:prelim}
To aid in presenting the proposed adaptive filtering approach, some preliminaries will be introduced in this section that briefly touch on the topics of hyperbolic conservation laws, minimum entropy principles, discontinuous spectral element methods, and modal filtering.

\subsection{Hyperbolic Systems and Entropy Principles}
The present work pertains to approximations of hyperbolic conservation laws (or hyperbolic components of mixed hyperbolic-parabolic conservation laws) of the form
\begin{equation}
    \label{eq:hyp_eq}
    \begin{cases}
    \partial_t \mathbf{u} + \boldsymbol{\nabla}{\cdot}{\mathbf{F}(\mathbf{u})} = 0, \\
    \mathbf{u}(\mathbf{x}, 0) = \mathbf{u}_0(\mathbf{x}), \quad \quad \mathrm{for}\  \mathbf{x} \in \Omega,
    \end{cases}
\end{equation}
where $\mathbf{u} \in \mathbb R^m$ is a vector-valued solution, $m$ is an arbitrary number of field variables, $\mathbf{F}(\mathbf{u}) \in \mathbb R ^{m \times d}$ is the flux, $\Omega \in \mathbb R^d$ is the domain, and $d$ is an arbitrary number of spatial dimensions. For simplicity, the domain is assumed to be periodic or the solution is compactly supported. If a solution of \cref{eq:hyp_eq} is an entropy solution \citep{Dafermos2010_45}, an entropy inequality of the form 
\begin{equation}\label{eq:entropy_law}
    \partial_t \sigma (\mathbf{u}) + \boldsymbol{\nabla}{\cdot}{\mathbf{\Sigma}(\mathbf{u})} \geqslant 0,
\end{equation} 
can be posed, for which the inequality holds given any entropy-flux pair $(\sigma, \mathbf{\Sigma})$ \citep{Harten1983} that satisfies 
    \begin{equation*}
        \pxvar{\mathbf{\Sigma }}{\mathbf{u}} = \pxvar{\sigma }{\mathbf{u}} \pxvar{\mathbf{F} }{\mathbf{u}}.
    \end{equation*} 
This notion of the entropy $\sigma$ is generally referred to as a \textit{numerical} entropy \citep{Tadmor1986}. We utilize a formulation of \cref{eq:entropy_law} that is of opposite sign in comparison to the work of \citet{Tadmor1986} for consistency with a physical entropy. In smooth regions, the inequality is satisfied exactly (i.e., the entropy functional obeys a conservation law), whereas in the vicinity of shocks, the left-hand side attains a strictly positive value (i.e., an entropy source exists). 

For certain systems, entropy solutions of \cref{eq:hyp_eq} satisfy a minimum principle on the entropy, i.e.,
\begin{equation}
    \sigma\left(\mathbf{u} (\mathbf{x}, t + \Delta t )\right) \geqslant \underset{\mathbf{x} \in \Omega}{\mathrm{min}}\ \sigma\left(\mathbf{u} (\mathbf{x}, t)\right),
\end{equation}
for all $\Delta t > 0$ \citep{Dafermos2010_128, Tadmor1986}. Assuming a finite propagation speed in hyperbolic systems, it is possible to restrict this condition over a local domain of influence to form a local minimum principle on the entropy.  For a given point $\mathbf{x}_0$, let $D_0$ be some local domain of influence over the interval $[t, t+ \Delta t]$ -- e.g., a closed $d$-ball centered on $\mathbf{x}_0$ with a radius of $\lambda_{\mathrm{max}} \Delta t$, where $\lambda_{\mathrm{max}}$ is some upper bound on the local maximum propagation speed of the system. A local minimum entropy principle can then be given as 
\begin{equation}
    \sigma\left(\mathbf{u} (\mathbf{x}_0, t + \Delta t )\right) \geqslant \underset{\mathbf{x} \in D_0}{\mathrm{min}}\ \sigma\left(\mathbf{u} (\mathbf{x}, t)\right).
\end{equation}
For numerical approximations of hyperbolic conservation laws, it can be advantageous to find schemes that can enforce this condition in scenarios where enforcing \cref{eq:entropy_law} is not feasible as it still guarantees some notion of physicality to the solution and can help alleviate instabilities in the numerical scheme \citep{Khobalatte1994}.

\subsection{Discontinuous Spectral Element Methods}

\begin{figure}[tbhp]
    \centering
    \adjustbox{width=0.3\linewidth,valign=b}{\input{./figs/fr_elems.tex}}
    \caption{\label{fig:fr_elems} Diagram of the interior solution points (red circles, $\mathbf{u}_i$), interior interface flux/solution points (red circles, blue outline, $\mathbf{u}_i^+$), and exterior interface flux points (blue circles, $\mathbf{u}_i^-$) for a triangular $\mathbb P_3$ element $\Omega_k$.  }
\end{figure}
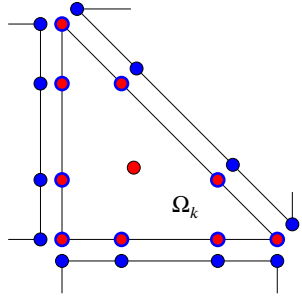

For a nodal discontinuous spectral element approximation of \cref{eq:hyp_eq} (e.g., discontinuous Galerkin \citep{Hesthaven2008DG}, flux reconstruction \citep{Huynh2007}, etc.), the domain $\Omega$ is partitioned into $N_e$ elements $\Omega_k$ such that $\Omega = \bigcup_{N_e}\Omega_k$ and $\Omega_i\cap\Omega_j=\emptyset$ for $i\neq j$, as shown in \cref{fig:fr_elems}. The approximate solution $\mathbf{u} (\mathbf{x})$ within each element $\Omega_k$ is given via a nodal approximation of the form
\begin{equation}
    \mathbf{u} (\mathbf{x}) = \sum_{i\in S} \mathbf{u}_i \phi_i(\mathbf{x}),
\end{equation}
where  $\mathbf{x}_i\ \forall \ i \in S$ is a set of solution nodes,  $\phi_i(\mathbf{x})$ are their associated nodal basis functions that possess the property $\phi_i(\mathbf{x}_j) = \delta_{ij}$, and $S$ is the set of nodal indices for the stencil. The shorthand notation $\mathbf{u}_i$ is used to denote the discrete nodal value of the solution (i.e., $\mathbf{u}_i = \mathbf{u}(\mathbf{x}_i)$). The order of the approximation, represented by $\mathbb P_p$ for some order $p$, is defined as the maximal order of $\mathbf{u} (\mathbf{x})$. 
The interaction between elements is communicated via the element interfaces $\partial \Omega$ using a set of interface nodes $\mathbf{x}_i \in \partial \Omega \ \forall \ i \in I$, where $I$ is a set of nodal indices for the interface stencil. In the work, we assume that the solution nodes are closed, such that the interface nodes are a subset of the solution nodes (i.e., $I \subset S$) to avoid issues regarding interpolation for discontinuous solutions. Modifications to allow for open solution nodes are possible but are outside of the scope of this work. Given an interface point $\mathbf{x}_i$ for some $i \in I$, let $\mathbf{u}_i^{-}$ denote the value of $\mathbf{u}(\mathbf{x})$ evaluated from the element of interest and let $\mathbf{u}_i^{+}$ denote the value of $\mathbf{u}(\mathbf{x})$ evaluated from the interface-adjacent element. 

For an arbitrary element $\Omega_k$ of some discontinuous spectral element approximation, the divergence of the flux at some point $\mathbf{x}_i$ can be approximated in terms of an interior flux component ($\boldsymbol{\nabla}{\cdot}{\mathbf{F}}_{\Omega_k}(\mathbf{u}_i)$) and an interface flux component ($\boldsymbol{\nabla}{\cdot}{\mathbf{F}}_{\partial \Omega_k}(\mathbf{u}_i)$). 
\begin{equation} 
    \boldsymbol{\nabla}{\cdot}{\mathbf{F}}(\mathbf{u}_i) \approx \boldsymbol{\nabla}{\cdot}{\mathbf{F}}_{\Omega_k}(\mathbf{u}_i) + \boldsymbol{\nabla}{\cdot}{\mathbf{F}}_{\partial \Omega_k}(\mathbf{u}_i).
\end{equation}
For the interior component, the divergence of the flux is calculated through a collocation projection of the flux onto the solution nodes, given as
\begin{equation}\label{eq:cij}
    \boldsymbol{\nabla}{\cdot}{\mathbf{F}}_{\Omega_k}(\mathbf{u}_i) = \sum_{j \in S} \mathbf{c}_{ij} \mathbf{F}(\mathbf{u}_j),
\end{equation}
for some discretization-dependent matrix $\mathbf{c}_{ij}$. For the interface component, the divergence of the flux is calculated as a function of both the interior and exterior interface values, given as
\begin{equation}\label{eq:cbarij}
    \boldsymbol{\nabla}{\cdot}{\mathbf{F}}_{\partial \Omega_k}(\mathbf{u}_i) = \sum_{j \in I} \overline{\mathbf{c}}_{ij} \overline{\mathbf{F}}(\mathbf{u}_j^{-}, \mathbf{u}_j^{+}, \mathbf{n}_j),
\end{equation}
where $\overline{\mathbf{c}}_{ij}$ is again some discretization-dependent matrix and $\overline{\mathbf{F}}(\mathbf{u}_j^{-}, \mathbf{u}_j^{+}, \mathbf{n}_j)$ is a common interface flux value dependent on the interior/exterior solution values ($\mathbf{u}_j^{-}$, $\mathbf{u}_j^{+}$) and their associated normal vector $\mathbf{n}_j$. This common interface flux is generally computed using exact or approximate Riemann solvers such as that of \citet{Rusanov1962} and \citet{Roe1981}. The semidiscretization of \cref{eq:hyp_eq} can then be given as
\begin{equation}
    \partial_t \mathbf{u}_i =  - \left( \boldsymbol{\nabla}{\cdot}{\mathbf{F}}_{\Omega_k}(\mathbf{u}_i) + \boldsymbol{\nabla}{\cdot}{\mathbf{F}}_{\partial \Omega_k}(\mathbf{u}_i) \right).
\end{equation}

We assume that the spatial scheme is chosen such that
\begin{equation}
    \partial_t \overline{\mathbf{u}} = - \int_{\partial \Omega_k}\overline{\mathbf{F}}\left(\mathbf{x} \right)\cdot\mathbf{n}(\mathbf{x})\ \mathrm{d}\mathbf{x} \approx - \sum_{j \in I} m_j \overline{\mathbf{F}}(\mathbf{u}_j^{-}, \mathbf{u}_j^{+}, \mathbf{n}_j)
\end{equation}
where $m_j$ is the corresponding quadrature weight for the point $\mathbf{x}_j$ and $\overline{\mathbf{u}}$ is the element-wise mean defined as
\begin{equation}
    \overline{\mathbf{u}} = \frac{1}{V_k}\int_{\Omega_k}\mathbf{u} (\mathbf{x})\ \mathrm{d}\mathbf{x} \quad \quad \mathrm{and} \quad \quad V_k = \int_{\Omega_k} \mathrm{d}\mathbf{x},
\end{equation}
for some arbitrary element $\Omega_k$. This relation is recovered for nodal discontinuous Galerkin approximations with appropriate quadrature and flux reconstruction schemes utilizing the equivalent discontinuous Galerkin correction functions \citep{Vincent2010}, the latter of which is the focus of this work. The ability of DSEM to preserve desirable properties on the element-wise mean is well documented in the literature due to its equivalency to first-order finite volume schemes \citep{Zhang2010, Zhang2011, Zhang2011b}. Without presenting a comprehensive proof, we assume that properties such as positivity of certain convex functionals of the solution and a discrete local minimum entropy principle are satisfied by the element-wise mean if one utilizes explicit strong stability preserving time integration under some Courant-Friedrichs–Lewy (CFL) condition \citep{Courant1928} with an appropriate choice for the Riemann solver. A more detailed description of these conditions and assumptions is presented in \cref{sec:implementation}.

\subsection{Modal Filtering}
The approximation of the solution given by a nodal basis can be equivalently expressed by a modal expansion as 
\begin{equation}
    \mathbf{u} (\mathbf{x})  = \sum_{i \in S} \widehat{\mathbf{u}}_i \psi_i(\mathbf{x}),
\end{equation}
where $\psi_i(\mathbf{x})\ \forall \ i \in S$ are a set of modal basis functions (e.g., Legendre polynomials, Koornwinder polynomials, etc.) and $\widehat{\mathbf{u}}_i$ are their corresponding modes. The modal basis is generally chosen such that the basis functions are orthogonal with respect to the inner product. With this formulation, a filtered solution can be defined as 
\begin{equation}
    \widetilde{\mathbf{u}}(\mathbf{x}) = \sum_{i \in S} H_i\left(\widehat{\mathbf{u}}_i\right) \psi_i(\mathbf{x}),
\end{equation}
where $H_i\left(\widehat{\mathbf{u}}_i\right)$ denotes some filtering operation applied to the modes. The filter function can be arbitrarily chosen, but must be dissipative, i.e.,
\begin{equation}
    \big | H_i\left(\widehat{\mathbf{u}}_i\right) \big | \leq \big | \widehat{\mathbf{u}}_i \big |,
\end{equation}
and conservative, i.e.,
\begin{equation}
    \frac{1}{V_k}\int_{\Omega_k}\widetilde{\mathbf{u}} (\mathbf{x})\ \mathrm{d}\mathbf{x}= \overline{\mathbf{u}}.
\end{equation}
For many filtering approaches, the filter tends to be more dissipative for higher frequency modes as spurious oscillations tend to manifest as high-frequency modes in the solution \citep{Persson2006}. 

%% file: figs/fr_elems.tex
     \begin{tikzpicture}[spy using outlines={rectangle, height=3cm,width=2.3cm, magnification=3, connect spies}]
		\begin{axis}[name=plot1,
		    axis line style={draw=none},
		    tick style={draw=none},
		    axis x line=left,
            axis y line=left,
            axis equal image,
            clip mode=individual,
    		xmin=-2,
    		xmax=1,
    		xticklabels={,,},
    		ymin=-1,
    		ymax=1,
    		yticklabels={,,},
    		style={font=\Large},
    		scale = 1]
    		
		    \draw[-] (axis cs:-1, -1) -- (axis cs:1, -1);
		    \draw[-] (axis cs:1, -1) -- (axis cs:-1, 1);
		    \draw[-] (axis cs:-1, 1) -- (axis cs:-1, -1);
		    
            \draw[-,fill=red] (axis cs:-0.333333333333333,-0.333333333333333) circle[radius=0.06];
            \draw[-,ultra thick, color=blue,fill=red] (axis cs:-1.0, 1.0) circle[radius=0.06];
            \draw[-,ultra thick,color=blue,fill=red] (axis cs:1.0, -1.0) circle[radius=0.06];
            \draw[-,ultra thick,color=blue,fill=red] (axis cs:-1.0, -1.0) circle[radius=0.06];
            \draw[-,ultra thick,color=blue,fill=red] (axis cs:0.447213595499958, -0.447213595499958) circle[radius=0.06];
            \draw[-,ultra thick,color=blue,fill=red] (axis cs:-1.0, -0.447213595499958) circle[radius=0.06];
            \draw[-,ultra thick,color=blue,fill=red] (axis cs:-0.447213595499958, -1.0) circle[radius=0.06];
            \draw[-,ultra thick,color=blue,fill=red] (axis cs:0.447213595499958, -1.0) circle[radius=0.06];
            \draw[-,ultra thick,color=blue,fill=red] (axis cs:-1.0, 0.447213595499958) circle[radius=0.06];
            \draw[-,ultra thick,color=blue,fill=red] (axis cs:-0.447213595499958, 0.447213595499958) circle[radius=0.06];

		    \draw[-] (axis cs:-1.2, 1) -- (axis cs:-1.5, 1);
		    \draw[-] (axis cs:-1.2, -1) -- (axis cs:-1.5, -1);
		    \draw[-] (axis cs:-0.86, 1.14) -- (axis cs:-0.36, 1.14);
		    
		    \draw[-] (axis cs:-1, -1.2) -- (axis cs:-1, -1.5);
		    \draw[-] (axis cs:1, -1.2) -- (axis cs:1, -1.5);
		    
		    \draw[-] (axis cs:1.14, -0.86) -- (axis cs:1.14, -0.56);
		    
		    \draw[-] (axis cs:-1, -1.2) -- (axis cs:1, -1.2);
		    \draw[-] (axis cs:1.14, -0.86) -- (axis cs:-0.86, 1.14);
		    \draw[-] (axis cs:-1.2, 1) -- (axis cs:-1.2, -1);
		    
            \draw[-,fill=blue] (axis cs:-1.2, 1.0) circle[radius=0.06];
            \draw[-,fill=blue] (axis cs:-1.2, -1.0) circle[radius=0.06];
            \draw[-,fill=blue] (axis cs:-1.2, -0.447213595499958) circle[radius=0.06];
            \draw[-,fill=blue] (axis cs:-1.2, 0.447213595499958) circle[radius=0.06];

            \draw[-,fill=blue] (axis cs:1.0, -1.2) circle[radius=0.06];
            \draw[-,fill=blue] (axis cs:-1.0, -1.2) circle[radius=0.06];
            \draw[-,fill=blue] (axis cs:-0.447213595499958, -1.2) circle[radius=0.06];
            \draw[-,fill=blue] (axis cs:0.447213595499958, -1.2) circle[radius=0.06];
            
            \draw[-,fill=blue] (axis cs:-0.86, 1.14) circle[radius=0.06];
            \draw[-,fill=blue] (axis cs:0.587213595499958, -0.307213595499958) circle[radius=0.06];
            \draw[-,fill=blue] (axis cs:-0.307213595499958, 0.587213595499958) circle[radius=0.06];
            \draw[-,fill=blue] (axis cs:1.14, -0.86) circle[radius=0.06];
            
            \node[align=left] at (0.15,-0.7) {$\Omega_k$};
            
		\end{axis}

	\end{tikzpicture}

%% file: methodology.tex
%!TEX root = ./main.tex
\section{Methodology}\label{sec:methodology}
To enforce certain desirable properties of the systems in question (at least in a discrete sense), it is possible to formulate them as convex constraints on the solution. A common constraint for physical systems is that some convex functional $\Gamma \left(\mathbf{u}\right)$ is non-negative across the domain (e.g., density and pressure for the Euler and Navier--Stokes equations, water height in the shallow water equations, etc.). This constraint is given by the condition 
\begin{equation}\label{eq:pp_con}
    \Gamma\left(\widetilde{\mathbf{u}}(\mathbf{x}_i) \right) \geq 0  \quad \forall \ i \in S,
\end{equation}
where the choice of functional(s) is dependent on the system in question. However, simply enforcing positivity of these functionals is usually not enough to ensure a well-behaved solution in the vicinity of a discontinuity. This constraint must generally be accompanied by a more restrictive condition, such as some local discrete minimum entropy principle on the solution, given by the condition
\begin{equation}\label{eq:mep_con}
    \sigma\left(\widetilde{\mathbf{u}}(\mathbf{x}_i)\right) \geq \sigma_{\min} \quad \forall \ i \in S,
\end{equation}
where $\sigma(\mathbf{u})$ is some convex entropy functional of the system in question and $\sigma_{\min}$ is some local minimum entropy to be defined in \cref{ssec:mep}. Enforcing this condition on the entropy tends to alleviate many of the issues regarding spurious oscillations in the vicinity of discontinuities \citep{Khobalatte1994}, although this is not always guaranteed if this condition is not enforced for \emph{every} possible entropy-flux pair.

\subsection{Adaptive Filtering}
To enforce these properties, we define an adaptive filtering operation with the goal of satisfying the constraints on the discrete filtered solution without the need for problem-dependent tunable parameters. The specific choice of filter is not particularly important as long as the filter meets the following criteria:
\begin{itemize}
    \item The filter is dependent on a single free parameter $\zeta$ (i.e., $H(\mathbf{u}) = H(\mathbf{u}, \zeta)$).  
    \item The filter is conservative (i.e., $\int_{\Omega_k}\widetilde{\mathbf{u}} (\mathbf{x})\ \mathrm{d}\mathbf{x} = \int_{\Omega_k}\mathbf{u} (\mathbf{x})\ \mathrm{d}\mathbf{x}$).
    \item There exists a minimum (or maximum) value of $\zeta$ such that the filter recovers the unfiltered solution (i.e., $\widetilde{\mathbf{u}}(\mathbf{x}) = \mathbf{u}(\mathbf{x})$).
    \item There exists a maximum (or minimum) value of $\zeta$ such that the filter recovers the mean mode (i.e., $\widetilde{\mathbf{u}}(\mathbf{x}) = \overline{\mathbf{u}}$).
\end{itemize}
The objective of the adaptive filter is to apply the minimum amount of filtering to the solution such that these constraints are met. We assume that for some arbitrary system, the discretization is chosen such that $\overline{\mathbf{u}}$ satisfies \cref{eq:pp_con,eq:mep_con}, an assumption that is explored in more detail in \cref{sec:implementation} for the specific systems in this work. From this assumption, it can be seen that there exists at least one value of $\zeta$ that recovers a filtered solution for which the constraints are satisfied exactly, i.e., \cref{eq:pp_con,eq:mep_con} are satisfied and at least one inequality becomes an equality. This value of $\zeta$ can be calculated through simple root-finding methods such as the bisection algorithm, and, in practice, is generally unique, although this is not guaranteed. 

In this work, a second-order exponential filter \citep{Glaubitz2017} is chosen, given by the filter function  
\begin{equation}
    H_i\left(\widehat{\mathbf{u}}_i \right) = \widehat{\mathbf{u}}_i\ e^{ -\zeta p_i^2},
\end{equation}
where $p_i$ is defined as the maximal order of the modal basis function $\psi_i(\mathbf{x})$. For this choice of filter, setting $\zeta = 0$ recovers the unfiltered solution and setting $\zeta = \infty$ recovers the mean mode, the latter of which can be approximated in a computational sense as $\zeta = \mathcal O \left(- \log ({\epsilon}) \right)$ for some value of machine precision $\epsilon$. With this formulation, we define the filter strength using the minimum value of $\zeta$ such that the solution abides by the positivity-preserving and discrete minimum entropy principle satisfying conditions, given by
\begin{equation}
    \quad \zeta = \underset{\zeta\ \geq\ 0}{\mathrm{arg\ min}} \ \ \mathrm{s.t.} \ \  \left [\Gamma\left(\widetilde{\mathbf{u}}(\mathbf{x}_i) \right) \geq 0, \ \sigma\left(\widetilde{\mathbf{u}}(\mathbf{x}_i)\right) \geq \sigma_{\min} \ \ \forall \ i \in S\right ].
\end{equation}
From a computational perspective, convergence to a local minima of $\zeta$ is sufficient in the case that there exist multiple values of $\zeta$ such that the constraints are satisfied exactly. A description of the approach for computing $\zeta$ as well as the choice of functionals to enforce constraints upon for the various systems is given in \cref{sec:implementation}. For sufficiently-resolved smooth solutions, the unfiltered solution is expected to already abide by these constraints \citep{Guermond2011}, and therefore the standard DSEM approximation would be recovered. The proposed filtering approach is hereafter referred to as entropy filtering as the method effectively filters the modes of the solution that contribute to the violation of a minimum entropy principle. 

\begin{remark}[Limiting] \label{rem:limit}
The proposed adaptive filtering operation can be considered somewhat similar to limiting-type approaches such as flux-corrected  transport \citep{Boris1997} and convex limiting \citep{Guermond2019}. In fact, for an unconventional filtering operation, given by
\begin{equation*}
    H_i\left(\widehat{\mathbf{u}}_i \right) =
    \begin{cases}
          \hphantom{\zeta}\widehat{\mathbf{u}}_i \quad \quad \mathrm{if\ } {p_i = 0}, \\
           \zeta \widehat{\mathbf{u}}_i\quad\quad \mathrm{else,}
    \end{cases}
\end{equation*} 
the filter recovers the approach of \citet{Zhang2010}, a linear convex limiting operation between a low-order solution ($\zeta = 0$) and a high-order solution ($\zeta = 1$). For a less trivial filter choice, this can be considered to be a nonlinear (and generally non-convex) limiting operation. 
\end{remark}

\subsection{Entropy Constraints}\label{ssec:mep}
Although positivity-preserving constraints are generally unequivocal for most physical systems, the notion of a minimum entropy principle is more ambiguous. Many systems are endowed with a multitude of entropy functionals $\sigma(\mathbf{u})$ \citep{Tadmor1986} and the choice of $\sigma_{\min}$ is not clearly defined. In this work, the minimum entropy principle is enforced on a numerical entropy that can be chosen arbitrarily from any entropy-flux pair ($\sigma, \boldsymbol{\Sigma}$) that satisfies \cref{eq:entropy_law} for the given system, the particular choice of which is posited to have a minor overall effect.

To calculate $\sigma_{\min}$, the discretization is assumed to be explicit in time under some standard CFL condition. With this assumption, the domain of influence of an arbitrary element $\Omega_k$ over a single temporal integration step can be considered to be strictly contained within the element and its direct Voronoi neighbors. Thus, $\sigma_{\min}$ can be defined using information only from an element and its direct neighbors. Let $\sigma^k_*$ be defined as the minimum entropy within an element $\Omega_k$, given as 
\begin{equation}
    \sigma^k_* = \underset{i \in S } {\min}\ \sigma\left({\mathbf{u}_k}(\mathbf{x}_i)\right),
\end{equation}
where $\mathbf{u}_k$ denotes the solution within the element $\Omega_k$. Furthermore, let $\mathcal A (k)$ be the set of element indices which are face-adjacent with $\Omega_k$, including $\Omega_k$ itself. The local minimum entropy $\sigma_{\min}^k$ associated with the element $\Omega_k$ is then calculated as
\begin{equation}
    \sigma_{\min}^k = \underset{i \in \mathcal A (k)} {\min}\ \sigma^i_*.
\end{equation}
For elements adjacent to boundaries, the entropy of the boundary state is used as the adjacent entropy value. Assuming that these minimum entropy values are calculated prior to a temporal integration step and the filter is applied afterwards, this formulation of the entropy constraint enforces that the local discrete minimum entropy is non-decreasing in time across its domain of influence.

%% file: implementation.tex
%!TEX root = ./main.tex
\section{Implementation}\label{sec:implementation}
\subsection{Governing Equations and Constraints}
The efficacy of the entropy filtering approach was evaluated on hyperbolic and mixed hyperbolic-parabolic conservation laws. For the hyperbolic system, the compressible Euler equations were chosen, given in the form of \cref{eq:hyp_eq} as 
\begin{equation}\label{eq:euler}
    \mathbf{u} = \begin{bmatrix}
            \rho \\ \boldsymbol{\rho v} \\ E
        \end{bmatrix} \quad  \mathrm{and} \quad \mathbf{F} = \begin{bmatrix}
            \boldsymbol{\rho v}\\
            \boldsymbol{\rho v}\otimes\mathbf{v} + P\mathbf{I}\\
        (E+P)\mathbf{v}
    \end{bmatrix},
\end{equation}
where $\rho$ is the density, $\boldsymbol{\rho v}$ is the momentum, $E$ is the total energy, $P = (\gamma-1)\left(E - \shalf\rho\mathbf{v}{\cdot}\mathbf{v}\right)$ is the pressure, and $\gamma = 1.4$ is the ratio of specific heat capacities for air. The symbol $\mathbf{I}$ denotes the identity matrix in $\mathbb{R}^{d\times d}$ and $\mathbf{v} = \boldsymbol{\rho v}/\rho$ denotes the velocity. Positivity constraints were placed on the density and pressure,
\begin{equation*}
    \Gamma_1 (\mathbf{u}) = \rho \quad \quad \mathrm{and} \quad \quad 
    \Gamma_2 (\mathbf{u}) = P,
\end{equation*} 
and the entropy functional was chosen as
\begin{equation*}
    \sigma = \rho \log(P\rho^{-\gamma}),
\end{equation*}
taken from the entropy-flux pair $(\sigma, \mathbf{v}\sigma)$.
Assuming an explicit strong stability preserving temporal integration scheme under some standard CFL condition with a solution that initially satisfies these constraints discretely, the element-wise mean at the next temporal step will satisfy these constraints if the interface fluxes are computed using an entropy-stable positivity-preserving Riemann solver \citep{Zhang2011b, Chen2017} (e.g., Godunov methods \citep{Godunov1959}, local Lax-Friedrichs flux \citep{Lax1954}, HLLC \citep{Toro1994} with appropriate wavespeed estimates). These properties of the element-wise mean were shown for discontinuous Galerkin approximations in the works of \citet{Zhang2010}, \citet{Zhang2011}, \citet{Zhang2011b}, and \citet{Chen2017}.

\begin{remark}[Source terms] \label{rem:source}
The extension of the proposed approach to hyperbolic systems with source terms is possible for discretizations that preserve the constraints on the element-wise mean. For positivity constraints, \citet{Zhang2011c} showed these properties under a potentially more restrictive time step condition. For the entropy constraint, the contribution of the source term to the entropy over the temporal integration step would have to be evaluated to augment the $\sigma_{\min}$ value, or if the entropy source is strictly positive, one may forego this modification and apply the proposed filter approach at the expense of a more-relaxed entropy constraint. 
\end{remark}

For the mixed hyperbolic-parabolic system, the compressible Navier--Stokes equations were chosen. As the assumption on the entropy of the element-wise mean is not necessarily satisfied with the mixed hyperbolic-parabolic discretization, the hyperbolic and parabolic components of the conservation law were isolated and treated separately using an explicit operator splitting approach \citep{Demkowicz1990}. For this system, the conservation law can be represented as
\begin{equation}\label{eq:hyppar_eq}
    \partial_t \mathbf{u}^n + \boldsymbol{\nabla}{\cdot}\left({\mathbf{F}_I(\mathbf{u}^n)} + {\mathbf{F}_V(\mathbf{u}^n)}\right) = 0,
\end{equation}
where the solution and the inviscid (hyperbolic) and viscous (parabolic) components of the flux, denoted by the subscripts $I$ and $V$, respectively, are given as 
\begin{equation}\label{eq:navierstokes}
    \mathbf{u} = \begin{bmatrix}
            \rho \\ \boldsymbol{\rho v} \\ E
        \end{bmatrix}, \quad  \mathbf{F}_I = \begin{bmatrix}
            \boldsymbol{\rho v}\\
            \boldsymbol{\rho v}\otimes\mathbf{v} + P\mathbf{I}\\
        (E+P)\mathbf{v}
    \end{bmatrix}, \quad \mathrm{and} \quad \mathbf{F}_V = \begin{bmatrix}
            0\\
            -\mu \left(\nabla \mathbf{v} + \nabla \mathbf{v}^T \right) \\
        -\mu \left(\nabla \mathbf{v} + \nabla \mathbf{v}^T \right)\mathbf{v} - \mu \frac{\gamma}{Pr} \nabla e 
    \end{bmatrix},
\end{equation}
where $\mu$ is the dynamic viscosity, $Pr = 0.73$ is the Prandtl number for air, and $e = \rho^{-1} (E - \shalf\rho \mathbf{v}{\cdot}\mathbf{v})$ is the specific internal energy. Identical positivity and entropy constraints are used for the Navier--Stokes equations as for the Euler equations. To enforce an entropy constraint for this system, at each temporal integration stage, the hyperbolic step is computed first, after which the entropy filter can be applied as with purely hyperbolic systems. The parabolic component of the temporal update is then added to the filtered solution. A final check is performed to ensure that the mixed hyperbolic-parabolic solution retains the positivity-preserving properties of the hyperbolic step. In the rare occasion that it does not, the filter is applied again using only positivity constraints. For both components of the flux, the boundary conditions were modified to ensure consistency with the system being solved (e.g., no-slip boundary conditions for the parabolic step were replaced with slip boundary conditions for the hyperbolic step).

\begin{remark}[Navier--Stokes equations] \label{rem:mixed}
The notion of a minimum entropy principle is satisfied by the Navier--Stokes equations when considering the thermodynamic entropy (see \citet{Tadmor1986}, Section 3). If the discretization for the parabolic component can be formed such that the minimum entropy principle is satisfied on the element-wise mean, the operator splitting approach can be neglected and the filter can be applied on the full hyperbolic-parabolic step with a significant reduction in computational cost. Alternatively, one may neglect the operator splitting approach without modifying the parabolic discretization to reduce the computational cost at the expense of the entropy constraint not necessarily being satisfied. 
\end{remark}

\subsection{Discretization and Computational Framework}
The proposed entropy filtering approach was implemented within PyFR \citep{Witherden2014}, a high-order unstructured flux reconstruction \citep{Huynh2007} (FR) solver that can target multiple compute architectures including CPUs and GPUs. The FR framework can be cast in the form of \cref{eq:cij} and \cref{eq:cbarij} as 
\begin{equation}
    \mathbf{c}_{ij} = 
    \begin{cases}
    
    \nabla \ell_j(\mathbf{x}_i), \hphantom{ - \nabla \mathbf{g}_j(\mathbf{x}_i)} \quad \quad \mathrm{if} \ j \notin I, \\
    \nabla \ell_j(\mathbf{x}_i) - \nabla \mathbf{g}_j(\mathbf{x}_i) \quad \quad \mathrm{else},
    \end{cases}
\end{equation}
and
\begin{equation}
    \overline{\mathbf{c}}_{ij} =  \nabla \mathbf{g}_j(\mathbf{x}_i),
\end{equation}
where $\ell_j(\mathbf{x})$ corresponds to the Lagrange interpolating polynomial for the $j$-th solution point and $\mathbf{g}_j(\mathbf{x})$ corresponds to the correction function for the $j$-th interface point. The correction functions \citep{Castonguay2011, Trojak2021} posses the properties that 
\begin{equation*}
    \mathbf{n}_i{\cdot} \mathbf{g}_j(\mathbf{x}_i) = \delta_{ij} \quad \mathrm{and} \quad \sum_{i \in I} \mathbf{g}_i (\mathbf{x}) \in \mathrm{RT}_{p},
\end{equation*}
where $\mathrm{RT}_{p}$ is the Raviart--Thomas space \citep{Raviart1977} of order $p$. In this work, the correction functions are chosen such as to recover the nodal discontinuous Galerkin approach \citep{Huynh2007, Hesthaven2008DG}. For a more in-depth overview and details on the extension to second-order PDEs, the reader is referred to \citet{Witherden2016} and the references therein. 

The solution nodes were distributed along the Gauss--Legendre--Lobatto quadrature points for tensor-product elements and the $\alpha$-optimized points \citep{Hesthaven2008DG} for simplex elements. Common interface flux values were computed using the HLLC Riemann solver \citep{Toro1994} for the inviscid fluxes and the BR2 approach \citep{Bassi2000} for the viscous fluxes. Temporal integration was performed using a three-stage third-order strong stability preserving (SSP) Runge--Kutta scheme \citep{Gottlieb2001}. For the modal basis, orthogonal polynomials with respect to the unit measure were used (i.e., Legendre basis for tensor-product elements, Proriol-Koornwinder-Dubiner-Owens basis for triangles, etc.).

\subsection{Filter Implementation}
 The implementation of the entropy filter was formulated as an element-wise scalar optimization problem. At each substage of the temporal integration method, a filtering operation was performed on the solution to enforce the positivity-preserving and minimum entropy constraints, the latter of which was computed using the solution at the previous substage. If the unfiltered solution satisfied the constraints, no filter was applied. Otherwise, the minimum necessary filter strength $\zeta$ was calculated via 20 iterations of a bisection approach. Faster convergence could be obtained using more sophisticated root bracketing methods such as the Brent or Illinois methods, but these approaches were not explored in this work. For the density and pressure constraints, a minimum value of $\rho_{\min} = P_{\min} = 10^{-8}$ was enforced to ensure a non-vacuum state for the Riemann solver, such that the positivity constraints were instead implemented as 
\begin{equation*}
    \Gamma_1 (\mathbf{u}) = \rho - \rho_{\min} \quad \quad \mathrm{and} \quad \quad 
    \Gamma_2 (\mathbf{u}) = P - P_{\min}.
\end{equation*} 
For the entropy constraint, a numerical tolerance of $\epsilon_\sigma = 10^{-4}$ was given, such that the constraint was instead implemented as 
\begin{equation}
    \sigma\left(\widetilde{\mathbf{u}}(\mathbf{x}_i)\right) \geq \sigma_{\min} - \epsilon_\sigma \quad \forall \ i \in S.
\end{equation}
This comparatively larger tolerance resulted in a slightly relaxed entropy constraint which was found to be beneficial for two reasons. First, the constraint was notably more prone to numerical precision issues due to the logarithm operation, particularly in the limit as $\rho \to \rho_{\min}$ or $P \to P_{\min}$ since small variations in these values could cause orders of magnitude more variation in the entropy. Secondly, marginally better resolution of flow features could be obtained by allowing slight undershoots in the entropy as strictly enforcing the entropy principle can degrade the accuracy of the solution \citep{Khobalatte1994, Guermond2019}. A more detailed description of the computational implementation including pseudo-code is presented in \cref{app:algo}.

%% file: results.tex
%!TEX root = ./main.tex
\section{Results}\label{sec:results}
The proposed entropy filtering approach was evaluated on a series of numerical experiments for the Euler and Navier--Stokes equations within a high-order flux reconstruction framework. For brevity, the solution for these systems is expressed in terms of a vector of primitive variables as $\mathbf{q}=[\rho,\mathbf{v},P]^T$.

\subsection{Euler Equations}
\subsubsection{Sod Shock Tube}
For an initial evaluation of the shock capturing capabilities of the proposed approach in the context of the Euler equations, the canonical case of the Sod shock tube was considered \citep{Sod1978}. The problem assesses the ability of the approach in resolving the three main features of the Riemann problem: shock waves, rarefaction waves, and contact discontinuities. The domain is set to $\Omega=[0,1]$ and the initial conditions are given as 

    \begin{equation*}
        \mathbf{q}(x,0) = \begin{cases}
            \mathbf{q}_l, &\mbox{if } x\leqslant 0.5,\\
            \mathbf{q}_r, &\mbox{else},
        \end{cases} \quad \mathrm{given} \quad \mathbf{q}_l = \begin{bmatrix}
            1 \\ 0 \\ 1
        \end{bmatrix}, \quad \mathbf{q}_r = \begin{bmatrix}
            0.125 \\ 0 \\ 0.1
        \end{bmatrix}.
    \end{equation*}

\begin{figure}[tbhp]
        \centering
        
        \subfloat[$\mathbb P_3$]{\label{fig:sod_p3}         
        \adjustbox{width=0.48\linewidth, valign=b}{\input{figs/sod_p3}}}
        ~
        \subfloat[$\mathbb P_5$]{\label{fig:sod_p5}         
        \adjustbox{width=0.48\linewidth, valign=b}{\input{figs/sod_p5}}}
        ~
        \newline
        \caption{\label{fig:sod} Density profile of the Sod shock tube problem at $t = 0.2$ computed using a $\mathbb P_3$ (left) and $\mathbb P_5$ (right) FR approximation with ${\sim}200$ degrees of freedom.}
    \end{figure}
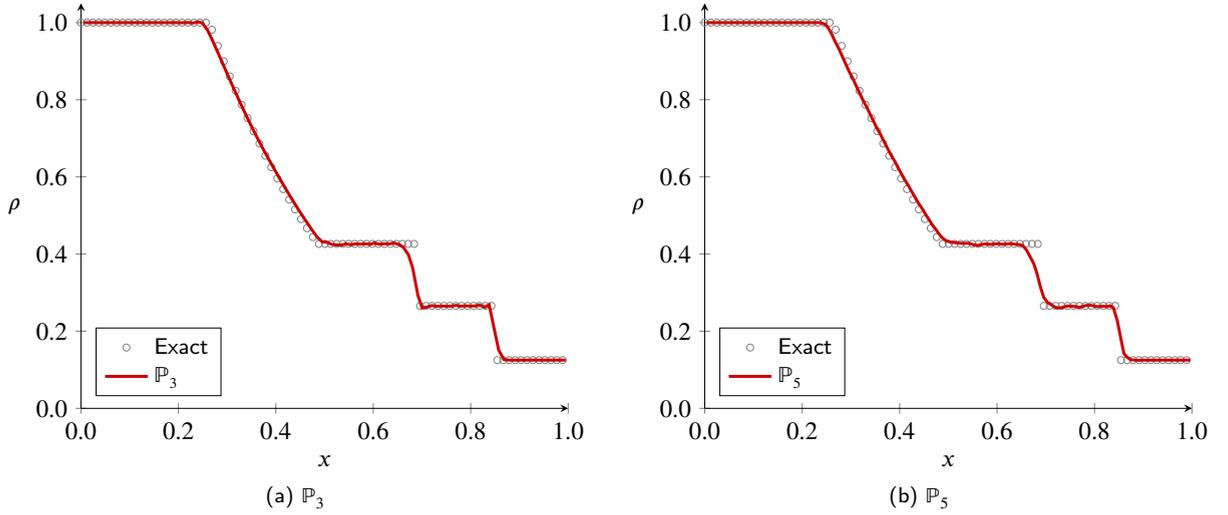

The density profiles at $t = 0.2$ as predicted by a $\mathbb P_3$ and $\mathbb P_5$ FR approximation with 200 degrees of freedom are shown in \cref{fig:sod}. For both approximation orders, the results showed good agreement with the exact solution, with excellent resolution of the rarefaction wave and shock wave and minimal dissipation around the contact discontinuity. Furthermore, negligible spurious oscillations were observed in the vicinity of discontinuities. For a fixed number of degrees of freedom, marginally better results were obtained using the lower-order $\mathbb P_3$ approximation than the higher-order $\mathbb P_5$ approximation, particularly around the contact discontinuity. This effect can be attributed to the proportionally higher number of elements available for lower-order approximations, giving a more localized approach for the filter. 

For a quantitative evaluation of the entropy filtering approach for discontinuous solutions, the convergence rates of the error against the exact solution were evaluated. For a given number of degrees of freedom $M$, the point-mean $L^1$ and $L^2$ norm of the density error was defined as 
    \begin{equation}
        \epsilon_{\rho_1} = \frac{1}{M}\sum_{i=0}^{M-1} \left|\rho (x_i) - \rho_{\mathrm{exact}} (x_i) \right| \quad \mathrm{and} \quad \epsilon_{\rho_2} = \sqrt{\frac{1}{M}\sum_{i=0}^{M-1} \left(\rho (x_i) - \rho_{\mathrm{exact}} (x_i) \right)^2},
    \end{equation}
respectively. The convergence rates of the density error with respect to the number of elements $N$ for various approximation orders are shown in \cref{tab:sod_error1} and \cref{tab:sod_error2}. The expected first-order convergence rate was generally obtained in both the $L^1$ and $L^2$ norm for all approximation orders, and the trend for the error was to decrease with increasing approximation order for a fixed number of elements. For a fixed number of degrees of freedom, shown in \cref{tab:sod_doferror1} and \cref{tab:sod_doferror2}, marginally lower error was generally obtained with a lower-order approximation due to the previously mentioned effects.

    \begin{figure}[tbhp]
    
        \centering
        \begin{tabular}{r | c c c c c c}\toprule
	        $N$ & $\mathbb{P}_2$ &$\mathbb{P}_3$ &$\mathbb{P}_4$ &$\mathbb{P}_5$ &$\mathbb{P}_6$ &$\mathbb{P}_7$ \\ \midrule
	        
            $40$ & \num{9.80e-03} & \num{8.57e-03} & \num{7.09e-03} & \num{8.08e-03} & \num{6.81e-03} & \num{6.96e-03} \\
            $80$ & \num{4.81e-03} & \num{4.30e-03} & \num{3.57e-03} & \num{4.15e-03} & \num{3.50e-03} & \num{3.61e-03} \\
            $160$ & \num{2.51e-03} & \num{2.33e-03} & \num{1.84e-03} & \num{2.16e-03} & \num{1.82e-03} & \num{1.99e-03} \\
            $320$ & \num{1.44e-03} & \num{1.30e-03} & \num{1.03e-03} & \num{1.18e-03} & \num{1.03e-03} & \num{1.14e-03} \\
            $640$ & \num{7.67e-04} & \num{6.03e-04} & \num{5.31e-04} & \num{6.32e-04} & \num{6.21e-04} & \num{6.38e-04} \\
            
            \midrule
        \textbf{RoC} & $0.91$ & $0.94$ & $0.93$ & $0.92$ & $0.87$ & $0.86$\\ 
        \end{tabular}
        \captionof{table}{\label{tab:sod_error1} Convergence of the $L^1$ norm of the density error with respect to mesh resolution $N$ for the Sod shock tube problem at $t = 0.2$ with varying orders. Rate of convergence shown on bottom.}
    \end{figure}
    
    \begin{figure}[tbhp]
    
        \centering
        \begin{tabular}{r | c c c c c c}\toprule
	        $N$ & $\mathbb{P}_2$ &$\mathbb{P}_3$ &$\mathbb{P}_4$ &$\mathbb{P}_5$ &$\mathbb{P}_6$ &$\mathbb{P}_7$ \\ \midrule
	        
            $20$ & \num{3.34e-03} & \num{2.04e-03} & \num{2.09e-03} & \num{1.58e-03} & \num{1.45e-03} & \num{1.39e-03} \\
            $40$ & \num{1.84e-03} & \num{1.08e-03} & \num{1.15e-03} & \num{8.12e-04} & \num{8.81e-04} & \num{7.04e-04} \\
            $80$ & \num{8.02e-04} & \num{4.39e-04} & \num{5.06e-04} & \num{3.11e-04} & \num{3.39e-04} & \num{2.82e-04} \\
            $160$ & \num{3.28e-04} & \num{2.43e-04} & \num{2.30e-04} & \num{1.75e-04} & \num{1.87e-04} & \num{1.45e-04} \\
            $320$ & \num{2.05e-04} & \num{1.45e-04} & \num{1.39e-04} & \num{1.04e-04} & \num{1.06e-04} & \num{8.59e-05} \\
            
            \midrule
        \textbf{RoC} & $1.05$ & $0.98$ & $1.01$ & $1.01$ & $0.98$ & $1.03$\\ 
        \end{tabular}
        \captionof{table}{\label{tab:sod_error2} Convergence of the $L^2$ norm of the density error with respect to mesh resolution $N$ for the Sod shock tube problem at $t = 0.2$ with varying orders. Rate of convergence shown on bottom.}
    \end{figure}

        \begin{figure}[tbhp]
    
        \centering
        \begin{tabular}{r | c c c c c c}\toprule
	        $M$ & $\mathbb{P}_2$ &$\mathbb{P}_3$ &$\mathbb{P}_4$ &$\mathbb{P}_5$ &$\mathbb{P}_6$ &$\mathbb{P}_7$ \\ \midrule
	        
            $100$ & \num{6.87e-03} & \num{9.51e-03} & \num{1.07e-02} & \num{1.27e-02} & \num{1.31e-02} & \num{1.72e-02} \\
            $200$ & \num{3.99e-03} & \num{5.37e-03} & \num{5.60e-03} & \num{7.34e-03} & \num{7.50e-03} & \num{9.70e-03} \\
            $400$ & \num{2.08e-03} & \num{2.72e-03} & \num{3.04e-03} & \num{3.88e-03} & \num{4.00e-03} & \num{5.65e-03} \\
            $800$ & \num{1.23e-03} & \num{1.46e-03} & \num{1.46e-03} & \num{2.08e-03} & \num{2.37e-03} & \num{2.90e-03} \\
            $1600$ & \num{6.60e-04} & \num{7.24e-04} & \num{9.00e-04} & \num{1.15e-03} & \num{1.24e-03} & \num{1.53e-03} \\
            
            \midrule
        \textbf{RoC} & $0.85$ & $0.93$ & $0.91$ & $0.88$ & $0.85$ & $0.87$\\ 
        \end{tabular}
        \captionof{table}{\label{tab:sod_doferror1} Convergence of the $L^1$ norm of the density error with respect to degrees of freedom $M$ for the Sod shock tube problem at $t = 0.2$ with varying orders. Rate of convergence shown on bottom.}
    \end{figure}
    
    \begin{figure}[tbhp]
    
        \centering
        \begin{tabular}{r | c c c c c c}\toprule
	        $M$ & $\mathbb{P}_2$ &$\mathbb{P}_3$ &$\mathbb{P}_4$ &$\mathbb{P}_5$ &$\mathbb{P}_6$ &$\mathbb{P}_7$ \\ \midrule
	        
            $100$ & \num{1.21e-03} & \num{1.42e-03} & \num{1.95e-03} & \num{1.78e-03} & \num{1.94e-03} & \num{2.34e-03} \\
            $200$ & \num{6.62e-04} & \num{7.21e-04} & \num{7.89e-04} & \num{9.23e-04} & \num{1.07e-03} & \num{1.07e-03} \\
            $400$ & \num{2.72e-04} & \num{3.24e-04} & \num{4.26e-04} & \num{3.75e-04} & \num{4.76e-04} & \num{5.68e-04} \\
            $800$ & \num{1.81e-04} & \num{1.71e-04} & \num{1.58e-04} & \num{1.92e-04} & \num{2.78e-04} & \num{2.44e-04} \\
            $1600$ & \num{9.47e-05} & \num{7.39e-05} & \num{1.15e-04} & \num{1.09e-04} & \num{1.46e-04} & \num{1.17e-04} \\
            
            \midrule
        \textbf{RoC} & $0.92$ & $1.06$ & $1.05$ & $1.03$ & $0.95$ & $1.08$\\ 
        \end{tabular}
        \captionof{table}{\label{tab:sod_doferror2} Convergence of the $L^2$ norm of the density error with respect to degrees of freedom $M$ for the Sod shock tube problem at $t = 0.2$ with varying orders. Rate of convergence shown on bottom.}
    \end{figure}
\subsubsection{Shu-Osher Problem}
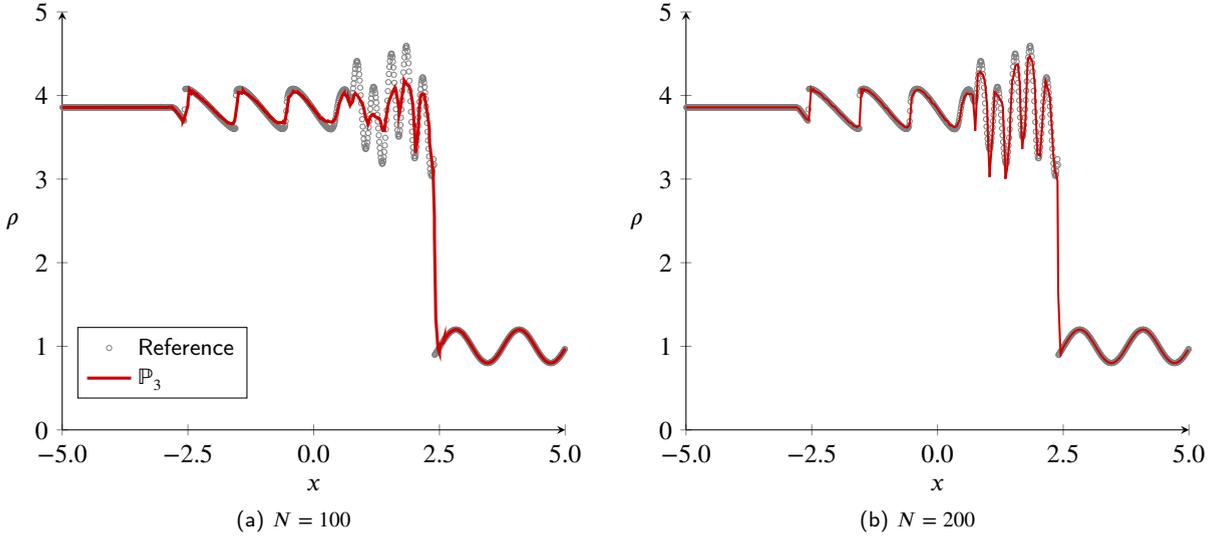
\begin{figure}[tbhp]
    \centering
    
    \subfloat[$N = 100$]{\label{fig:so_100}         
    \adjustbox{width=0.48\linewidth, valign=b}{\input{figs/shuosher100}}}
    ~
    \subfloat[$N = 200$]{\label{fig:so_200}         
    \adjustbox{width=0.48\linewidth, valign=b}{\input{figs/shuosher200}}}
    ~
    \newline
    \caption{\label{fig:shu} Density profile of the Shu-Osher problem at $t = 1.8$ computed using a $\mathbb P_3$ FR approximation with 100 (left) and 200 (right) elements.}
\end{figure}

To assess the effects of the entropy filter for more complex problems including shock waves and smooth oscillatory behavior, the case of \citet{Shu1988} was considered. The problem is solved on the domain $\Omega = [-5, 5]$ with the initial conditions
    \begin{equation*}
         \mathbf{q}(x,0) =  \begin{cases}
            \mathbf{q}_l, &\mbox{if } x\leqslant -4, \\
            \mathbf{q}_r, &\mbox{else},
        \end{cases} \quad \mathrm{given} \quad
        \mathbf{q}_l = \begin{bmatrix}
            3.857143 \\ 2.629369 \\ 10.333333
        \end{bmatrix}, \quad
        \mathbf{q}_r = \begin{bmatrix}
            1 + 0.2\sin{5x} \\ 0\\ 1
        \end{bmatrix}.
    \end{equation*}
The problem consists of a shock front propagating through a sinusoidally-perturbed density field, the interaction between which can induce instabilities in the flow field. However, these instabilities can be erroneously damped by overly dissipative shock capturing schemes. The predicted density profile at $t = 1.8$ computed using a $\mathbb P_3$ FR approximation with 100 and 200 elements is shown in \cref{fig:shu}. A reference solution was computed using a highly-resolved exact Godunov-type solver \citep{Toro1997_4}. The results show good resolution of the leading and trailing shock waves, and no spurious oscillations were observed. The instabilities in the field aft of the leading shock were also well-resolved, particularly with increasing resolution. 

\subsubsection{Isentropic Euler Vortex}
For an extension to two-dimensional problems, the entropy filtering approach was initially applied to a smooth solution to verify the accuracy of the underlying DSEM was not detrimentally affected. The isentropic Euler vortex problem \citep{Shu1998} was used as its analytic solution can be utilized to evaluate the convergence of the error. The initial conditions of the problem are given as

    \begin{equation*}
        \mathbf{q}(\mathbf{x}, 0)  = \begin{bmatrix}
            p^\frac{1}{\gamma} \\
            V_x + \frac{S}{2 \pi R} (y-y_0)\phi(r) \\[4pt]
            V_y - \frac{S}{2 \pi R} (x-x_0)\phi(r) \\[4pt]
            \frac{1}{\gamma M^2} \left(1 - \frac{S^2 M^2 (\gamma-1)} {8 \pi^2}\phi(r)^2\right)^\frac{\gamma}{\gamma-1}
        \end{bmatrix}, \quad \mathrm{where} \quad r = \|\mathbf{x}-\mathbf{x}_0\|_2 \quad \mathrm{and} \quad \phi(r) = \exp{\left(\frac{1-r^2}{2R^2}\right)},
    \end{equation*}
with the parameters $S = 13.5$ denoting the strength of the vortex, $R = 1.5$ the radius, $V_x = 0$, $V_y = 1$ the advection velocities, and $M = 0.4$ the freestream Mach number. The domain was set to $\Omega = [-10, 10]^2$, and a series of uniform quadrilateral meshes of size $N \times N$ was generated with periodic boundary conditions. After a single pass-through of the domain, the $L^2$ norm of the density error was calculated as
    \begin{equation}
        \epsilon_{\rho_2} = \sqrt{ \frac{1}{|\Omega|} \int_{\Omega} (\rho - \rho_{\mathrm{exact}})^2 \ \mathrm{d}{\mathbf{x}}},
    \end{equation}
with the quadrature calculated on $(2p)^2$ Gauss-Legendre nodes. This error is tabulated in \cref{tab:icv_error} for a series of experiments with varying mesh resolution and approximation order. The convergence rate of the error was generally in the range of $p$ to $p+1$ for the varying approximation orders, on par with the theoretical rate of $p+1$. These findings indicate that the entropy filter does not appreciably degrade the accuracy of the underlying DSEM for smooth solutions. 

    \begin{figure}[tbhp]
        \centering
        \begin{tabular}{r | c c c c c c }\toprule
	        $N$ & $\mathbb{P}_2$ &$\mathbb{P}_3$ &$\mathbb{P}_4$ &$\mathbb{P}_5$&$\mathbb{P}_6$&$\mathbb{P}_7$ \\ \midrule
	        
	        $20$ & - & - & - & - & \num{2.59E-05} & \num{5.37E-06}  \\ 
	        
	        $25$ & - & - & \num{7.80E-04} & \num{7.73E-05} & \num{6.84E-06} & \num{9.85E-07}     \\
	        
            $33$ & \num{1.80E-02}  & \num{1.79E-03} & \num{2.50E-04} & \num{1.11E-05} & \num{1.26E-06} & \num{9.58E-08}   \\
            
            $40$ & \num{1.10E-02} & \num{7.58E-04} & \num{1.08E-04} & \num{2.86E-06} & \num{3.75E-07} & \num{1.97E-08}    \\ 
            
            $50$ & \num{6.30E-03} & \num{3.02E-04} & \num{4.03E-05} & \num{7.52E-07}  & - & -   \\
            
            $67$ & \num{2.86E-03} & \num{1.05E-04} & - & -  & - & -   \\\midrule
        \textbf{RoC} & $2.59$& ${4.00}$& ${4.27}$& ${6.73}$ & ${6.10}$& ${8.13}$ \\ 
        \end{tabular}
        \captionof{table}{\label{tab:icv_error} Convergence of the $L^2$ norm of the density error with respect to mesh resolution $N$ for the isentropic Euler vortex problem with varying approximation order. Rate of convergence shown on bottom.}
    \end{figure}

To verify the computational efficiency of the entropy filter for smooth solutions where the unfiltered solution remains stable, the compute time for a solution without filtering and with a varying number of filter iterations was compared. The comparison was performed on an NVIDIA V100 GPU over 5 flow-throughs of the domain on an $N = 40$ mesh with a $\mathbb P_3$ approximation and $\Delta t = 2{\cdot}10^{-4}$. For 5, 10, and 20 iterations of the filter, the relative computational cost increase (as measured by relative time-to-solution) was 1.1\%, 1.7\%, and 2.0\%, respectively. These findings indicate that the entropy filter does not have a notable detrimental impact on the computational efficiency of the scheme for smooth solutions as it is primarily inactive. Furthermore, it also indicates that further improvements to the efficiency are possible through more sophisticated iterative solvers that require fewer iterations.
    
\subsubsection{Double Mach Reflection}
The double mach reflection problem of \citet{Woodward1984} was subsequently used to evaluate the ability of the entropy filter to resolve strong discontinuities in multiple dimensions. This case consists of a Mach 10 shock impinging on a $30$ degree ramp and results in multiple strong shock-shock and shock-contact interactions. The problem is solved on the domain $\Omega = [0,4] \times [0,1]$ with the initial conditions

    \begin{equation*}
         \mathbf{q}(\mathbf{x},0) =  \begin{cases}
            \mathbf{q}_l, &\mbox{if } x < 1/6 + \tan (30^\circ)y,\\
            \mathbf{q}_r, &\mbox{else},
        \end{cases} \quad \mathrm{given} \quad
        \mathbf{q}_l = \begin{bmatrix}
            8 \\ 7.14471\\ -4.125 \\ 116.5
        \end{bmatrix}, \quad
        \mathbf{q}_r = \begin{bmatrix}
            1.4 \\ 0\\ 0 \\ 1
        \end{bmatrix}.
    \end{equation*}
At the left boundary and the bottom boundary for $x < 1/6$, the solution was set to the post-shock state $\mathbf{q}_l$. No-slip adiabatic wall boundary conditions were applied for the bottom boundary for $x \geqslant \frac{1}{6}$. At the right boundary, the solution was set to the pre-shock state $\mathbf{q}_r$. For the top boundary, the exact solution is enforced, given as
    \begin{equation*}
        \mathbf{q}(\mathbf{x}, t)|_{y=1} = \begin{cases}
                \mathbf{q}_l, &\mbox{if } x\leqslant 1/6  + \tan (30^\circ)y + \frac{10}{\cos(30^\circ)}t,\\
                \mathbf{q}_r, &\mbox{else}.
            \end{cases}
    \end{equation*}

\begin{figure}[htbp!]
        \centering
        \adjustbox{width=0.72\linewidth,valign=b}{\includegraphics[width=\textwidth]{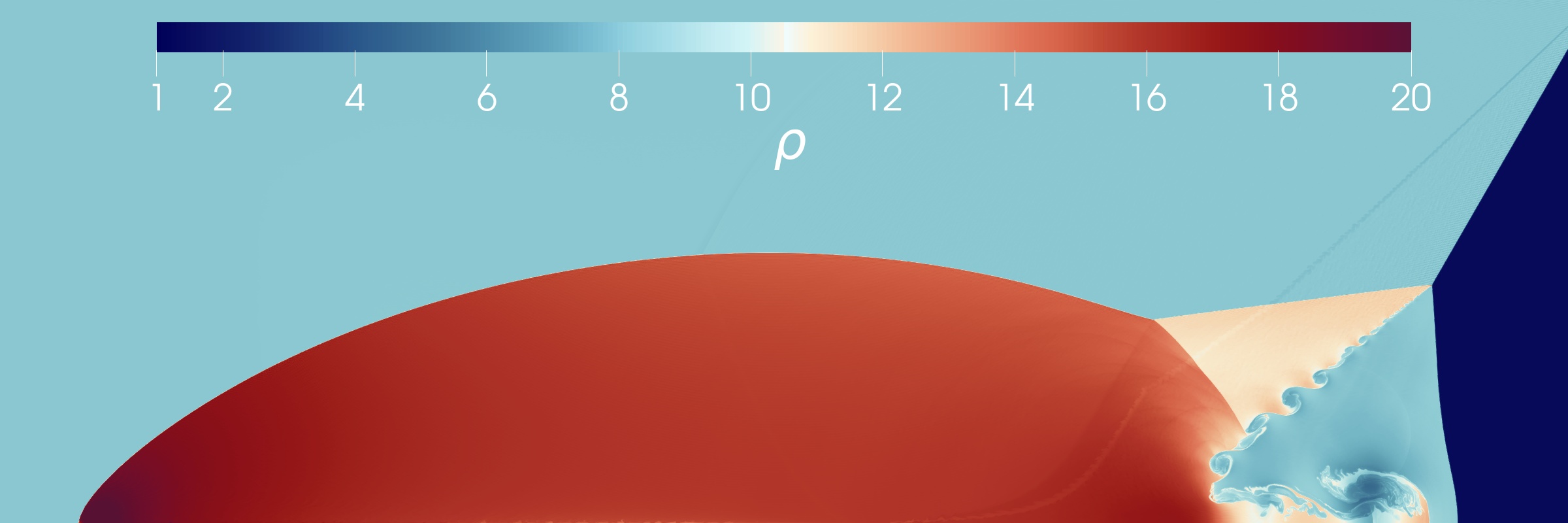}}
         \adjustbox{width=0.24\linewidth,valign=b}{\includegraphics[width=\textwidth]{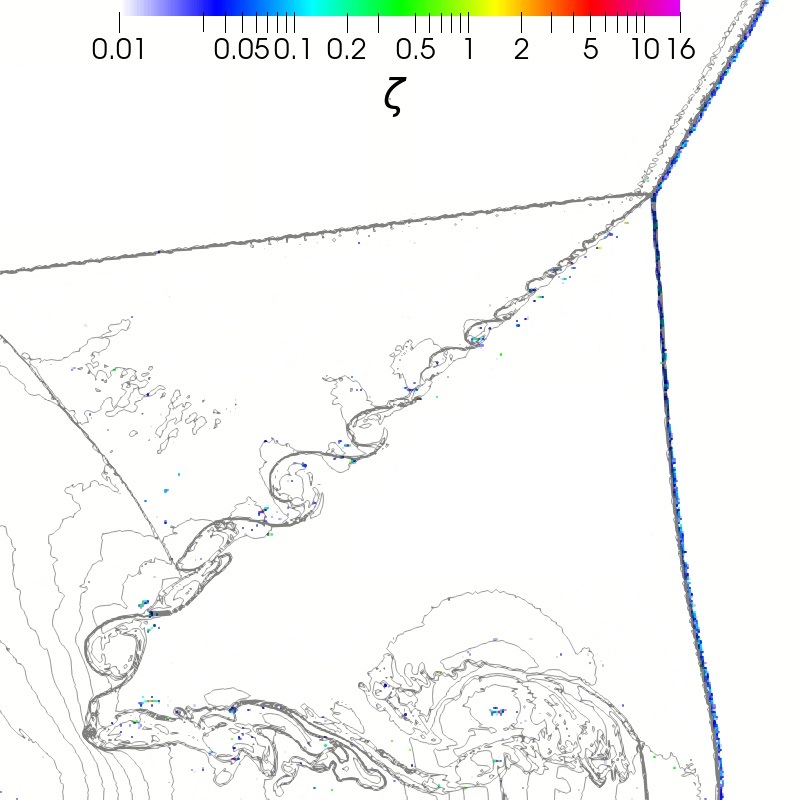}}
        \newline
        \caption{\label{fig:doublemach} (Left) Contours of density on the subregion $[0,3]\times[0,1]$ for the double Mach reflection problem at $t = 0.2$ using a $\mathbb P_3$ FR approximation with a $2400 \times 600$ mesh. (Right) Enlarged view of the distribution of the filter parameter $\zeta$ overlaid on isocontours of density.}
    \end{figure}

The contours of density at $t = 0.2$ as predicted by a $\mathbb P_3$ FR scheme on a $2400 \times 600$ mesh are shown in \cref{fig:doublemach}. The results show sub-element resolution of discontinuities without the presence of spurious oscillations. Furthermore, the application of the filter did not excessively dissipate the Kelvin-Helmholtz instabilities along the contact line, indicating that the filter is not erroneously dissipating small-scale features. To verify this, the distribution of the filter parameter $\zeta$ is also shown in \cref{fig:doublemach} overlaid on the isocontours of density. The filter was primarily active in the leading shock fronts, with minimal activation within the small-scale structures along contact line. Furthermore, even along the shock front, a relatively small value of $\zeta$ was required, significantly less than the value corresponding to the recovery of the mean mode. 

\subsubsection{Kelvin--Helmholtz Instability}
\begin{figure}[htbp!]
        \centering
        \subfloat[$N = 100^2$]{\adjustbox{width=0.4\linewidth,valign=b}{\includegraphics[width=\textwidth]{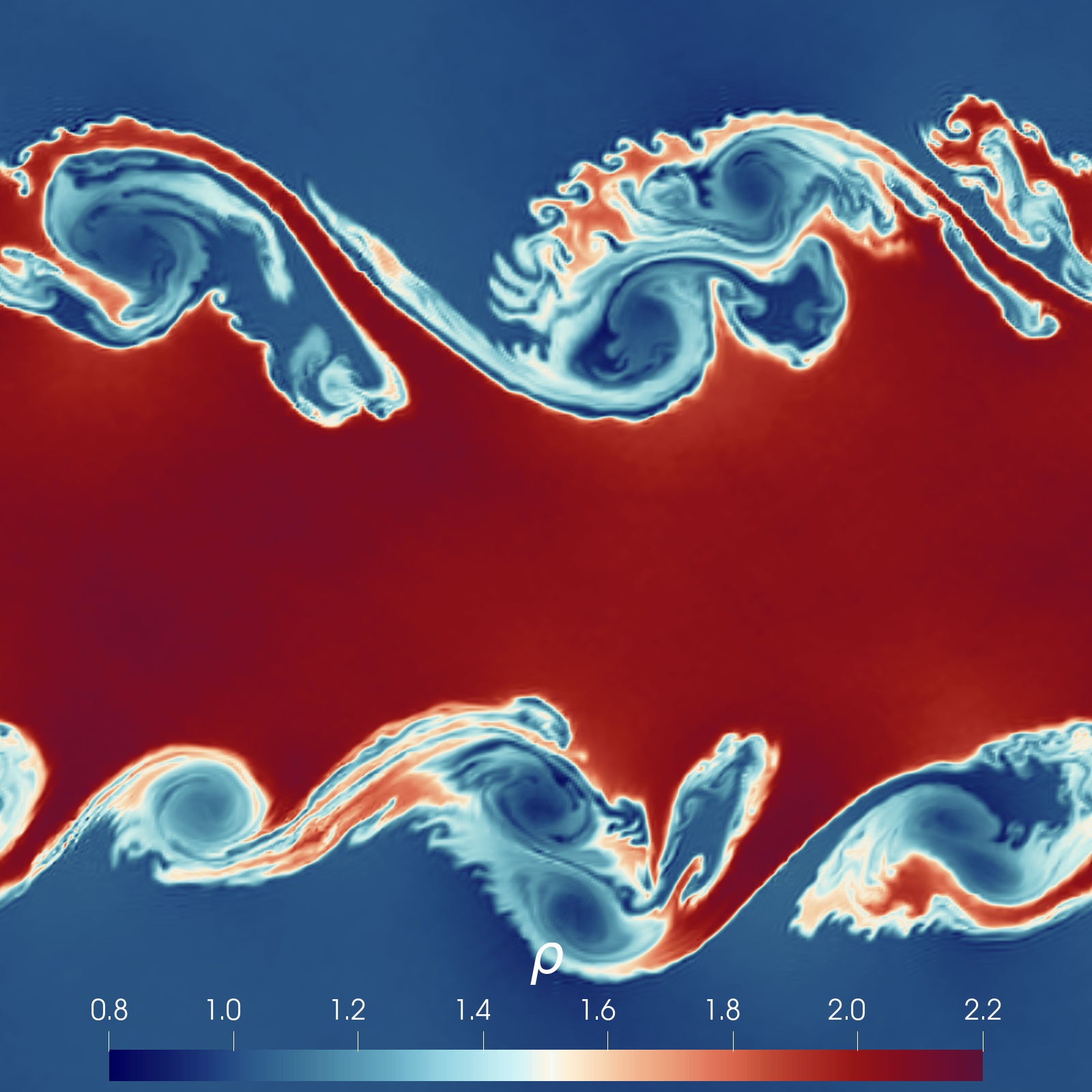}}}
        ~
        \subfloat[$N = 200^2$]{\adjustbox{width=0.4\linewidth,valign=b}{\includegraphics[width=\textwidth]{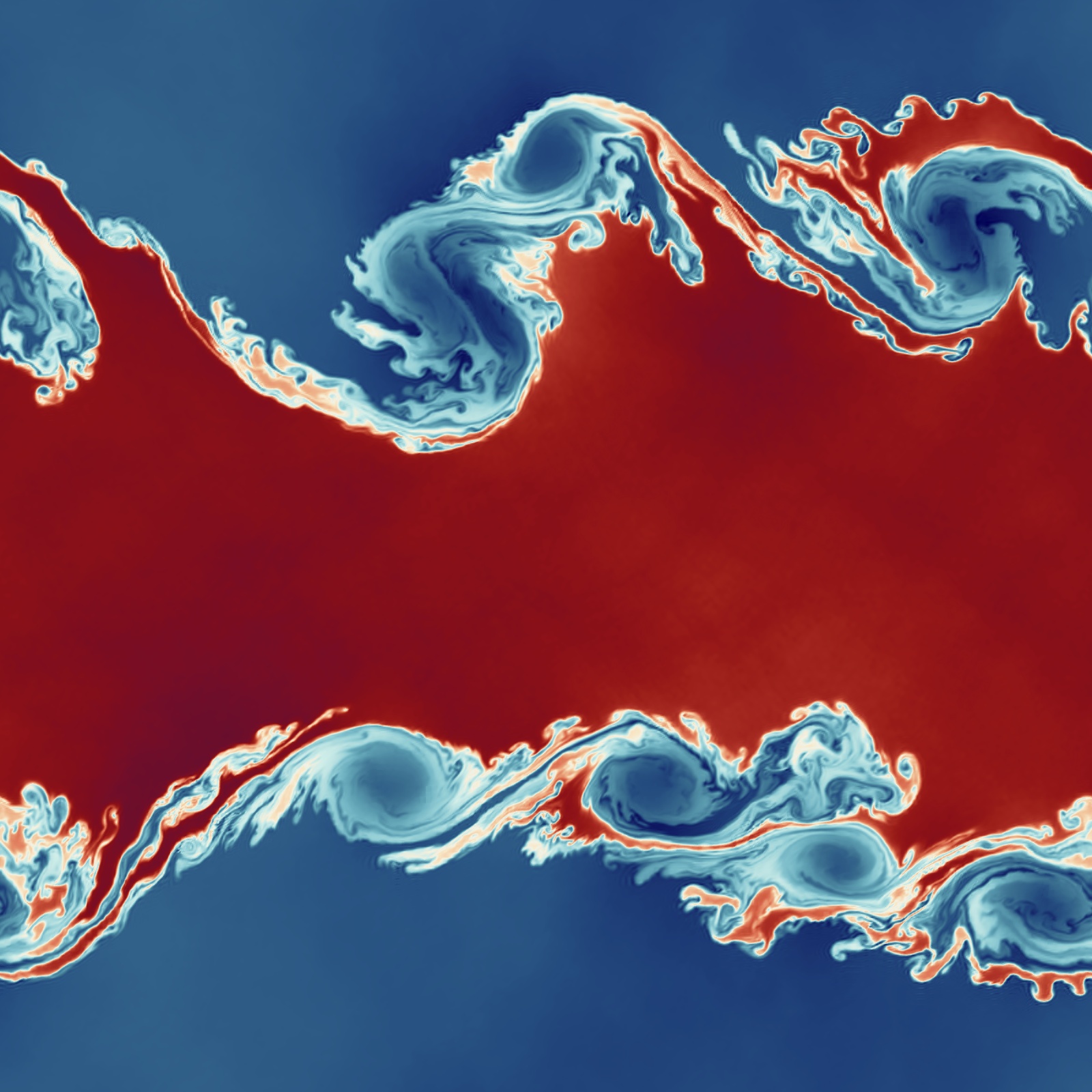}}}
        \newline
        \subfloat[$N = 400^2$]{\adjustbox{width=0.4\linewidth,valign=b}{\includegraphics[width=\textwidth]{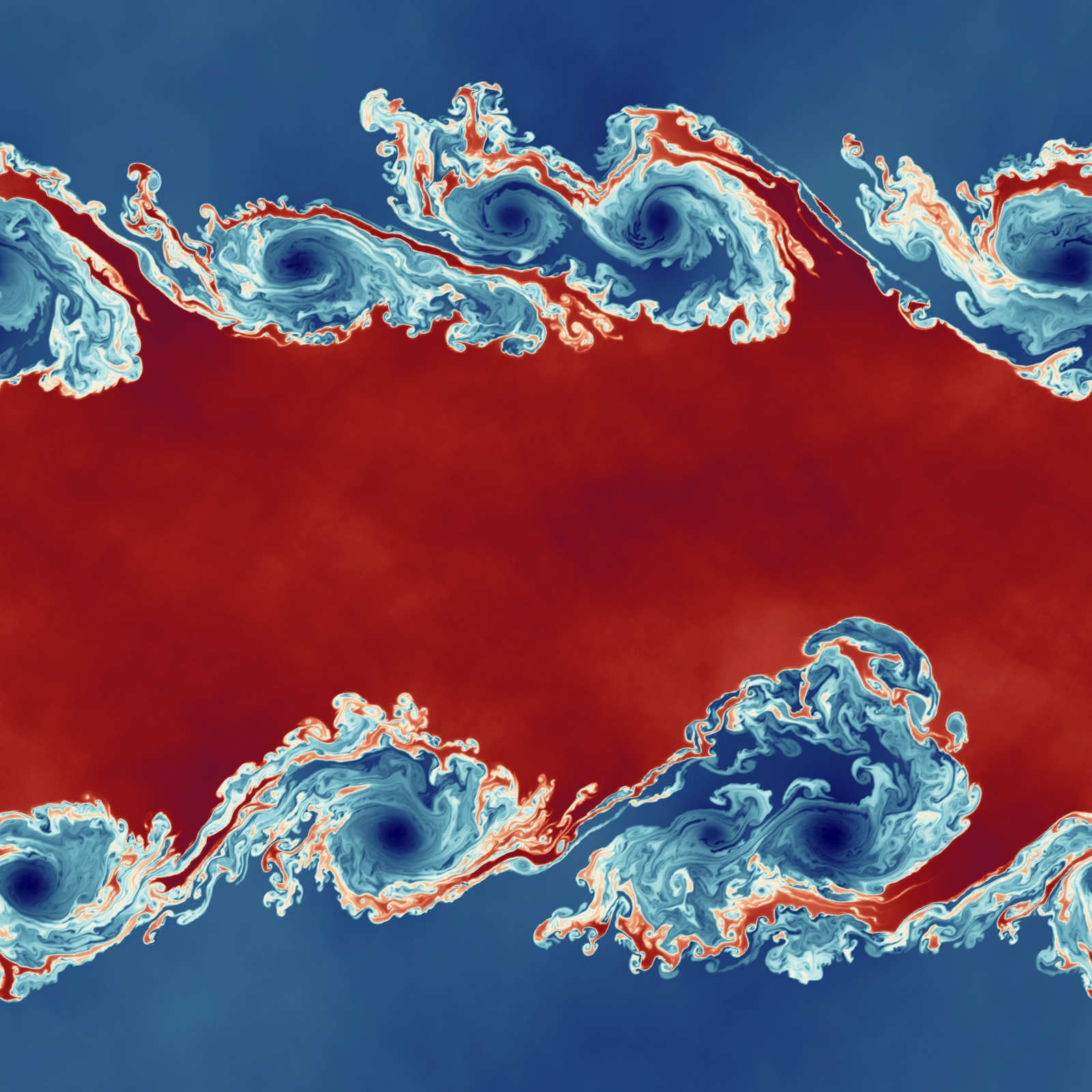}}}
        ~
        \subfloat[$N = 800^2$]{\adjustbox{width=0.4\linewidth,valign=b}{\includegraphics[width=\textwidth]{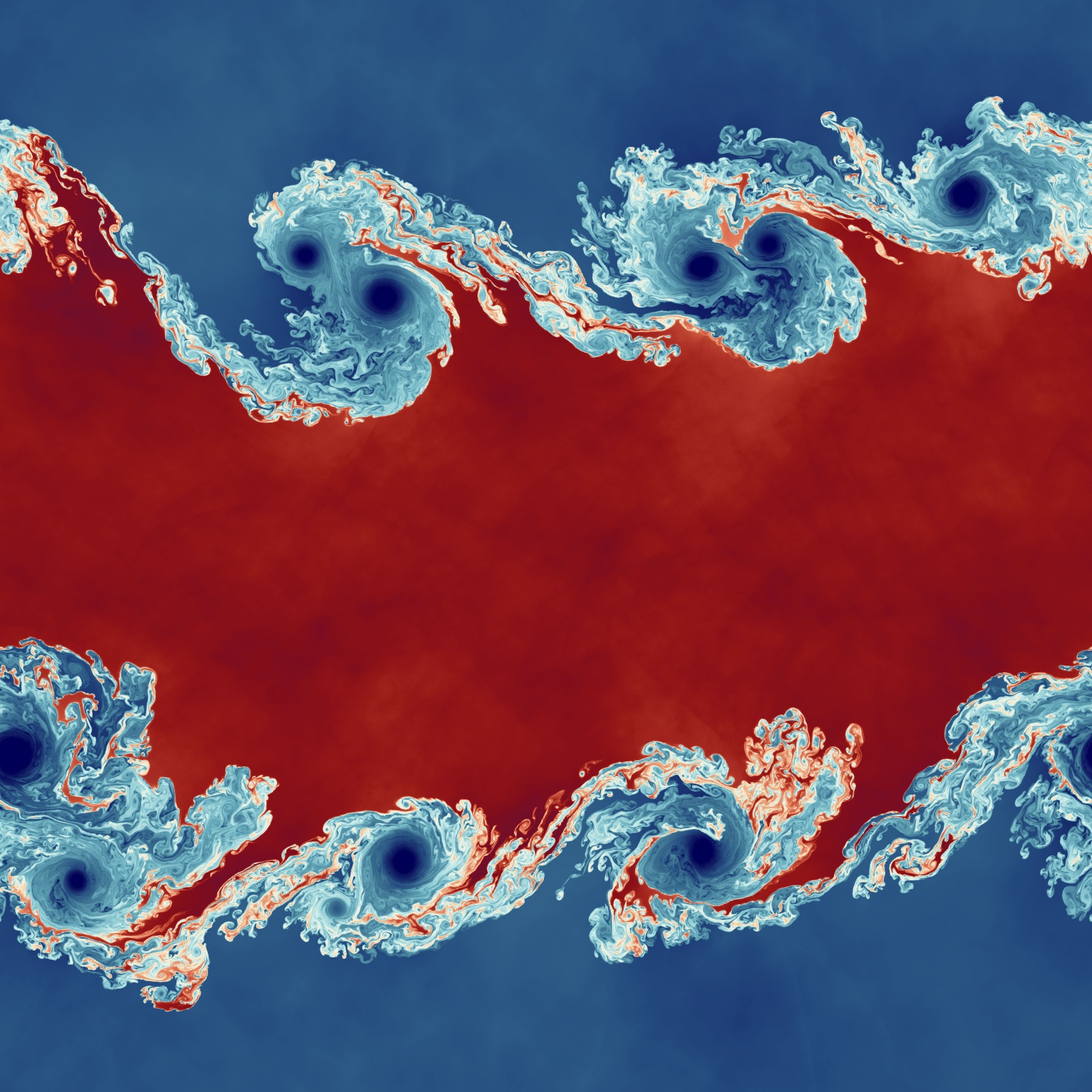}}}
        \newline
        \caption{\label{fig:kh} Contours of density for the Kelvin-Helmholtz instability problem at $t = 2$ using a $\mathbb P_4$ FR approximation with varying mesh resolutions.}
    \end{figure}

To further evaluate the ability of the entropy filtering approach to resolve small-scale flow features in the vicinity of discontinuities, the roll-up of a Kelvin-Helmholtz instability was simulated. The problem is solved on the domain $\Omega = [-0.5, 0.5]^2$ with the initial conditions 
    \begin{equation*}
         \mathbf{q}(\mathbf{x},0) =  \begin{cases}
            \mathbf{q}_l, &\mbox{if } |y| \leqslant 0.25,\\
            \mathbf{q}_r, &\mbox{else},
        \end{cases} \quad \mathrm{given} \quad
        \mathbf{q}_l = \begin{bmatrix}
            2 \\ 0.5\\ 0 \\ 2.5
        \end{bmatrix}, \quad
        \mathbf{q}_r = \begin{bmatrix}
            1 \\ -0.5\\ 0 \\ 2.5
        \end{bmatrix}.
    \end{equation*}
No explicit seeding of the instabilities was performed -- the instabilities originate from roundoff errors in the solver. The contours of density at $t = 2$ as predicted by a $\mathbb P_4$ FR scheme with various mesh resolutions are shown in \cref{fig:kh}. The roll-up of the vortices was well-resolved by the approach, with complex small-scale vortical structures beginning to appear with increasing mesh resolution. No spurious oscillations were observed in the vicinity of the discontinuities, and more subtle features such as pressure waves were not excessively dissipated. 

\subsubsection{Mach 800 Astrophysical Jet}
As a verification of the positivity-preserving properties of the filter for very extreme conditions, the case of high-speed astrophysical jets is considered. The test case, introduced by \citet{Balsara2012}, consists of a Mach 800 jet in the presence of an ambient gas. The problem setup is identical to the work of \citet{Wu2018} in which a provably-positive third-order discontinuous Galerkin approach is used, although the magnetic field is neglected in this case. A half domain $\Omega = [0, 0.5] \times [0, 1.5]$ is considered with symmetry (slip adiabatic) boundary conditions along the $y$-axis. The domain is filled with an ambient gas of density $0.1 \gamma$, zero velocity, and unit pressure. For the $y=0$ boundary, the inlet region is defined on $x \leqslant 0.05$, and the solution is set to $\mathbf{q} = [\gamma, 0, 800, 1]^T$ which yields a Mach number of 800 with respect to the inflow gas. The remaining boundary conditions are set as free. 

\begin{figure}[htbp!]
        \centering
        \subfloat[$N = 200 \times 600$]{\label{fig:jet_low} \adjustbox{width=0.48\linewidth,valign=b}{\includegraphics[width=\textwidth]{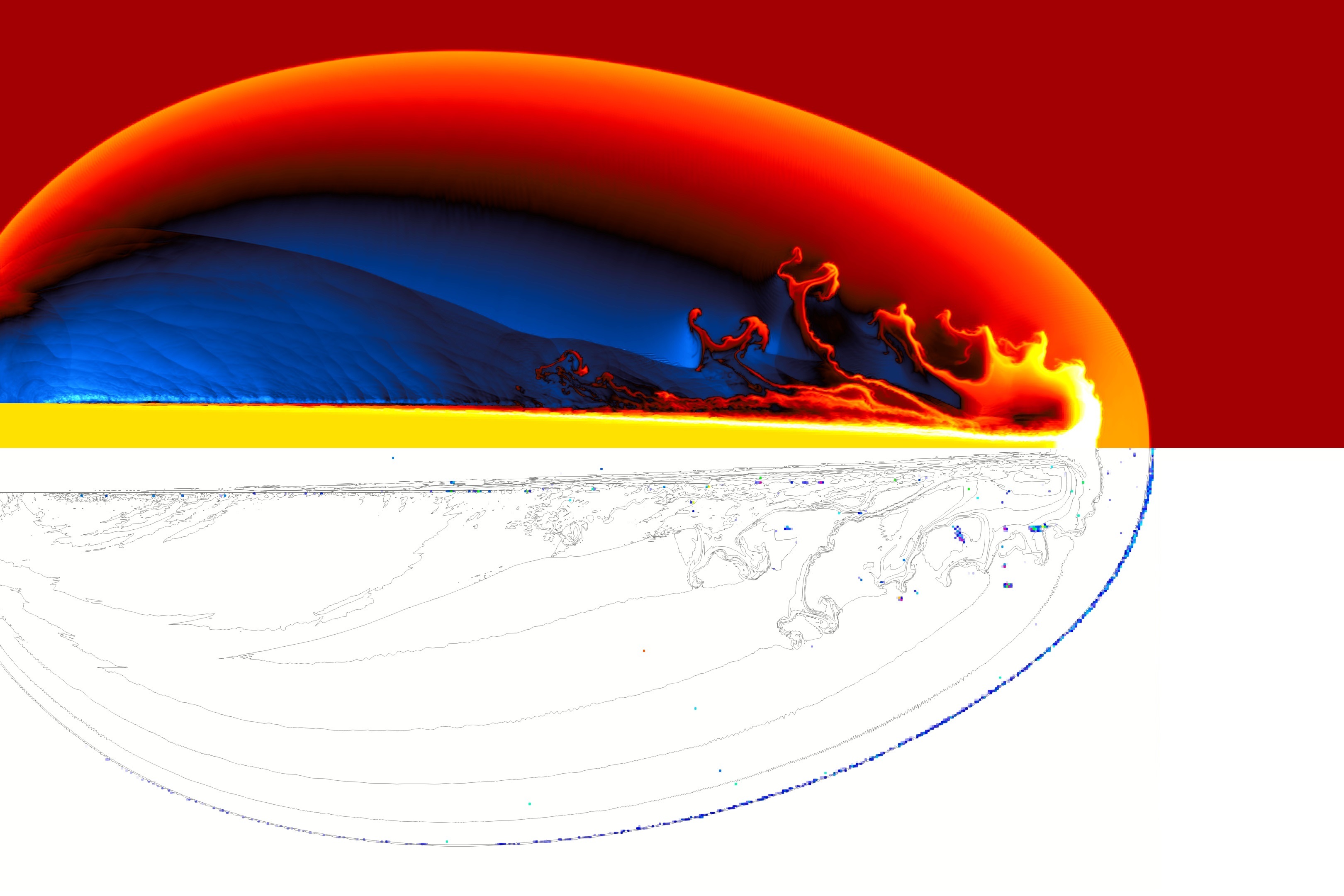}}}
        ~
        \subfloat[$N = 800 \times 2400$]{\label{fig:jet_high} \adjustbox{width=0.48\linewidth,valign=b}{\includegraphics[width=\textwidth]{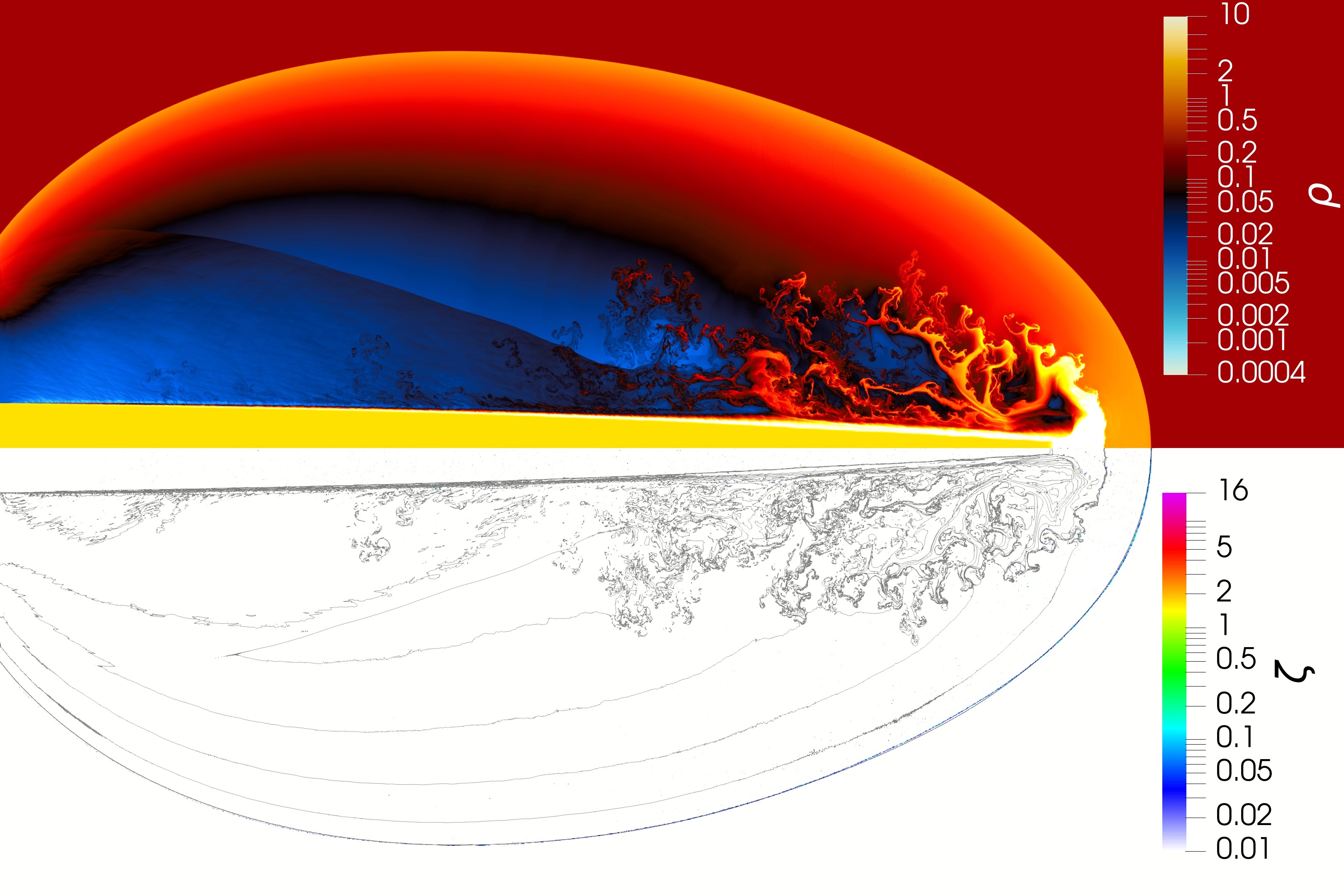}}}
        ~
        \newline
        \caption{\label{fig:jet} Contours of density (top) and distribution of the filter parameter $\zeta$ (bottom) for the Mach 800 astrophysical jet problem at $t = 0.002$ using a $\mathbb P_3$ FR approximation with a $200 \times 600$ mesh (left) and $800 \times 2400$ mesh (right). Contours are reflected about the $y$-axis. View is reoriented such that the $+y$ direction is shown left-to-right.}
    \end{figure}

The contours of density at $t = 0.002$ as predicted by a $\mathbb P_3$ FR scheme with a coarse ($200 \times 600$) and fine ($800 \times 2400$) mesh are shown in \cref{fig:jet}. Excellent resolution of the leading shock wave was obtained, and the small-scale structures in the vicinity of the cocoon/jet interface were well-resolved even with the coarse mesh. Additionally, the distribution of the filter parameter $\zeta$ is shown in \cref{fig:jet} overlaid on the isocontours of density. For both the coarse and fine mesh, the filter was primarily active at the leading shock region. Some regions of activation within the cocoon/jet interface were observed for the coarse mesh, but this behavior was reduced with increasing resolution, such that minimal activation away from the leading shock was observed with the fine mesh.

Given a similar case setup, a comparison can be made between the proposed approach on the coarse mesh and the mildly-magnetized results of \citet{Wu2018} (Fig. 6.i), which were obtained by the linear limiting approach of \citet{Zhang2010}. Significantly better resolution was obtained using the entropy filter, particularly with regards to the small-scale features near the cocoon/jet interface, although some of this may be attributed to the stabilizing effect of the magnetic field in their case. However, this suggests that the nonlinear limiting performed by the entropy filter offers noticeable advantages in comparison to a linear limiting approach with minimal computational overhead.

\subsection{Navier--Stokes Equations}
\subsubsection{Taylor--Green Vortex}
The entropy filtering approach was then extended to turbulent compressible flows through the Navier--Stokes equations. To verify that the filter does not unnecessarily dissipate small-scale turbulent fluctuations in the absence of shocks, the proposed approach was applied to the subsonic Taylor--Green vortex at a Reynolds number of 1600, a canonical fluid dynamics problem for studying vortex dynamics and turbulent transition and decay \citep{Taylor1937}. The problem is solved on the periodic domain $\Omega = [-\pi, \pi]^3$ with the initial conditions
    \begin{equation*}
         \mathbf{q}(\mathbf{x},0) =  \begin{bmatrix}
            1 \\ \hphantom{-}\sin(x)\cos(y)\cos(z)\\ -\cos(x)\sin(y)\cos(z) \\ 0\\ P_0 + \frac{1}{16}\left (\cos(2x) + \cos(2y) \right)\left (\cos(2z + 2) \right)
        \end{bmatrix},
    \end{equation*}
where $P_0 = 1/\gamma M^2$ for a reference Mach number $M = 0.08$. Given the unit density and velocity, the dynamic viscosity $\mu$ is set to $1/1600$ to recover a Reynolds number of 1600. 

The quantity of interest for this problem is the dissipation rate of the kinetic energy in the flow. The non-dimensional integrated kinetic energy can be defined as 
\begin{equation}
    K(t) = \frac{1}{V}\int_{\Omega} \frac{1}{2}\rho \mathbf{v}{\cdot} \mathbf{v} \ \mathrm{d}{\mathbf{x}},
\end{equation}
where $V = 8\pi^3$ is the volume of the domain. From this, a dissipation rate based on the kinetic energy can be calculated as
\begin{equation}
    \varepsilon_K = \frac{\mathrm{d} K}{\mathrm{d} t}.
\end{equation}
A similar measure of the dissipation can be obtained through a scaled form of the non-dimensional integrated enstrophy, defined as
\begin{equation}
    \varepsilon_E = \frac{ \beta}{V}\int_{\Omega} \frac{1}{2}\rho \boldsymbol{\omega}{\cdot} \boldsymbol{\omega} \ \mathrm{d}{\mathbf{x}},
\end{equation}
where $\boldsymbol{\omega}$ is the vorticity and $\beta = 2 \mu$ is the scaling factor. For purely incompressible flows, these two quantities are equal, but for compressible flows, they differ by the contribution of the deviatoric strain and pressure dilatation to the dissipation. At low Mach numbers, the enstrophy-based dissipation can reasonably approximate the kinetic energy-based dissipation for well-resolved flows.

\begin{figure}[tbhp]
    
        \subfloat[$\mathbb P_3$, $96^3$-$160^3$ DoF]{\label{fig:tgv_href}         
        \adjustbox{width=0.48\linewidth, valign=b}{\input{figs/tgv_dedt}}}
        ~
        \subfloat[120 DoF, $\mathbb P_3$-$\mathbb P_5$]{\label{fig:tgv_pref}         
        \adjustbox{width=0.48\linewidth, valign=b}{\input{figs/tgv_enst}}}
        ~
        \newline
        \caption{\label{fig:tgv} Dissipation measured by kinetic energy (black) and enstrophy (red) for the Taylor-Green vortex using a $\mathbb P_3$ FR approximation with varying DoF (left) and ${\sim}120$ DoF with varying approximation order (right). DNS results of \citet{vanRees2011} (obtained by private communication) shown for reference.}
    \end{figure}
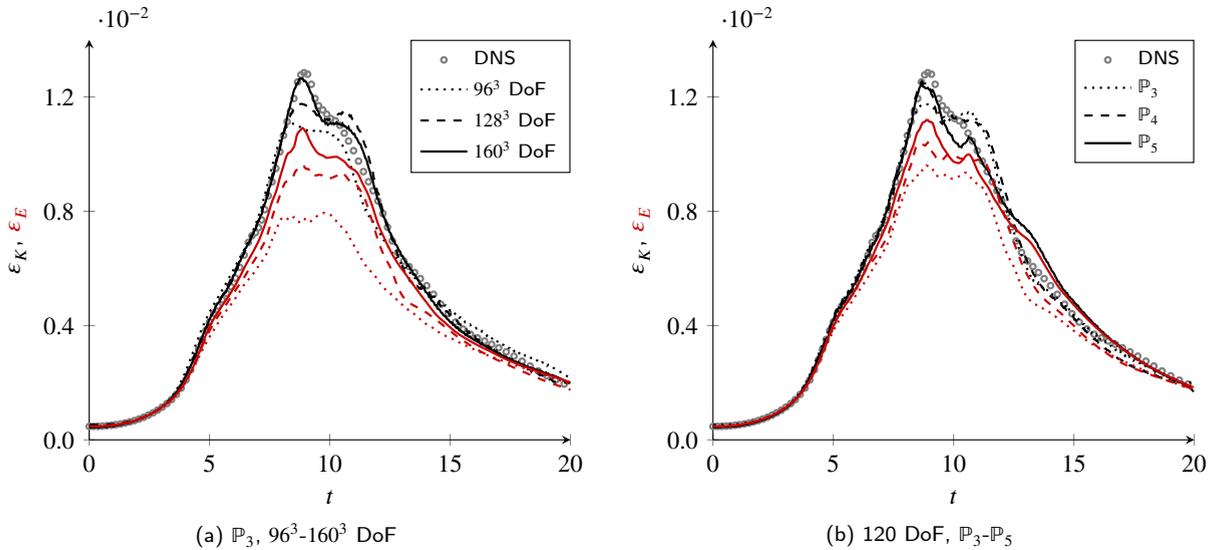

The prediction of these two quantities as computed by the entropy filtering approach with varying resolution and approximation order is shown in \cref{fig:tgv} in comparison to the DNS results of \citet{vanRees2011}. For visualization of the kinetic energy-based dissipation, a moving-average smoothing operation was performed prior to computing the temporal derivative to reduce oscillations. For a $\mathbb P_3$ FR approximation, the number of degrees of freedom was varied from $96^3$-$160^3$. Relatively good agreement was observed between the kinetic energy-based dissipation and the reference data across this entire range of resolution, with the most resolved case showing negligible deviation from the reference. For the enstrophy-based dissipation, the profiles evidently showed convergence to the reference with increasing resolution, and for a given resolution, showed similar results to unfiltered approaches \citep{Spiegel2016}. When fixing the degrees of freedom to ${\sim}120$ and varying the approximation order, improvements in the prediction of both the kinetic energy- and enstrophy-based dissipation were observed with increasing approximation order. Additionally, nearly identical results were obtained without filtering (see \citet{Trojak2021ACM}, Fig. 3b), which supports the presumption that the filter is predominantly inactive for turbulent flows without discontinuities.

\subsubsection{Viscous Shock Tube}
To assess the efficacy of the entropy filtering approach for predicting shock-boundary layer interactions, a two-dimensional viscous shock tube was simulated. The problem, introduced by \citet{Daru2000}, is solved on the half-domain $\Omega = [0,1]\times[0, 0.5]$ with the initial conditions
    \begin{equation*}
         \mathbf{q}(\mathbf{x},0) =  \begin{cases}
            \mathbf{q}_l, &\mbox{if } x \leqslant 0.5,\\
            \mathbf{q}_r, &\mbox{else},
        \end{cases} \quad \mathrm{given} \quad
        \mathbf{q}_l = \begin{bmatrix}
            120 \\ 0\\ 0 \\ 120/\gamma
        \end{bmatrix}, \quad
        \mathbf{q}_r = \begin{bmatrix}
            1.2 \\ 0\\ 0 \\ 1.2/\gamma
        \end{bmatrix}.
    \end{equation*}
For the top wall, a slip adiabatic wall boundary condition is applied to enforce symmetry. No slip adiabatic wall conditions are applied for the remaining walls. This case contains the standard features of the Riemann problem, namely a rarefaction wave, a contact discontinuity, and a shock wave, in addition to viscous interactions between these features and the wall. This interaction forms a complex lambda shock impinging on a viscous boundary layer, and as a result, presents a suitable test case for shock-boundary layer interactions. 

    \begin{figure}[htbp!]
        \centering
        \subfloat[$\mu = 1{\cdot}10^{-3}$]{\label{fig:st_1k} \adjustbox{width=0.48\linewidth,valign=b}{\includegraphics[width=\textwidth]{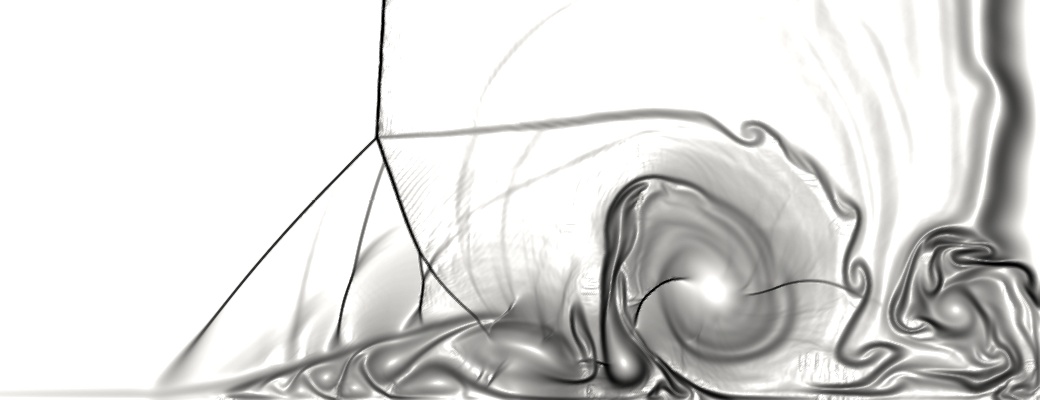}}}
        ~
        \subfloat[$\mu = 5{\cdot}10^{-4}$]{\label{fig:st_2k} \adjustbox{width=0.48\linewidth,valign=b}{\includegraphics[width=\textwidth]{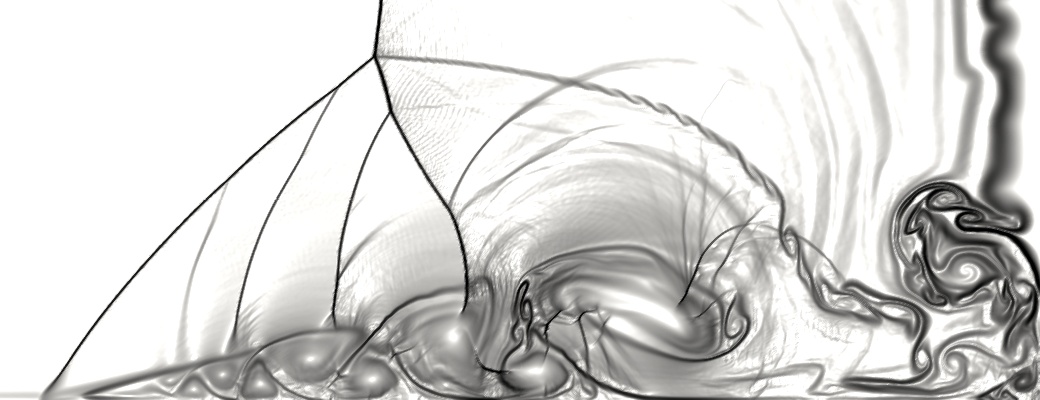}}}
        ~
        \newline
        \subfloat[$\mu = 2{\cdot}10^{-4}$]{\label{fig:st_5k} \adjustbox{width=0.48\linewidth,valign=b}{\includegraphics[width=\textwidth]{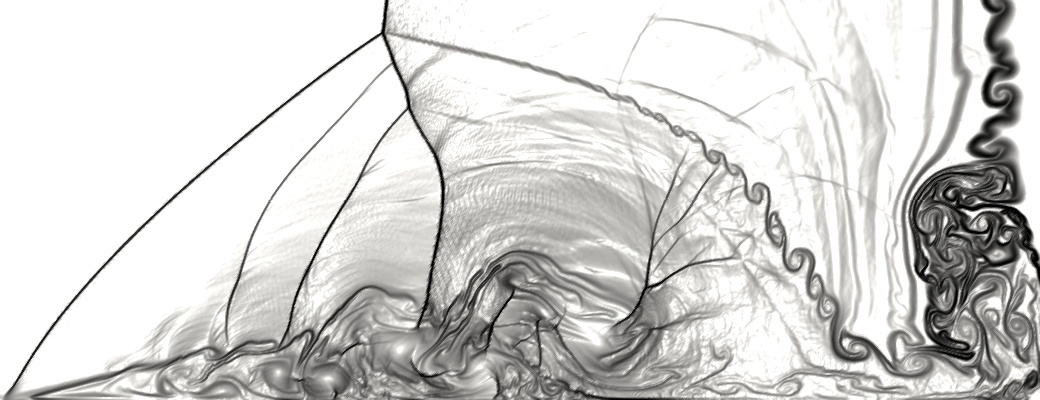}}}
        ~
        \subfloat[$\mu = 1{\cdot}10^{-4}$]{\label{fig:st_10k} \adjustbox{width=0.48\linewidth,valign=b}{\includegraphics[width=\textwidth]{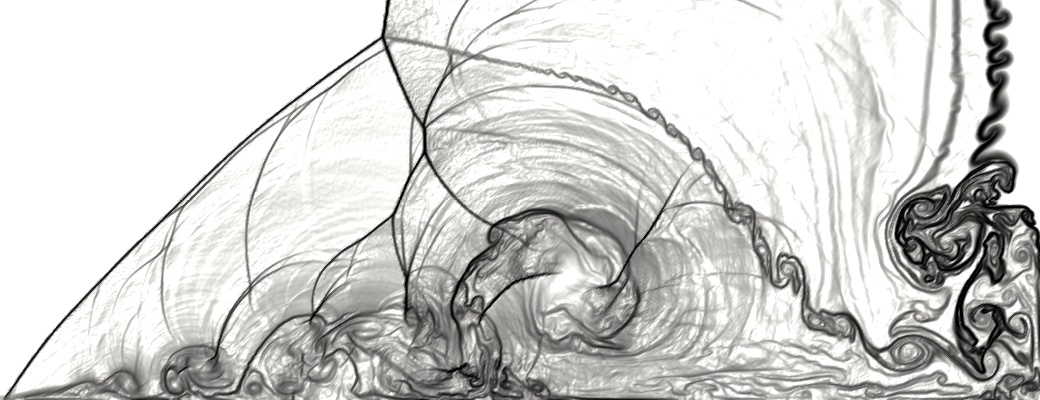}}}
        ~
        \newline
        \caption{\label{fig:shocktube_re} Schlieren-type representation of the density gradient norm for the viscous shock tube problem using a $\mathbb P_4$ FR approximation with a $800 \times 400$ mesh for varying values of the dynamic viscosity $\mu$. Contours shown at $t = 1$ on the subregion $[0.35,1]\times[0,0.25]$.}
    \end{figure}
    
Numerical Schlieren results at $t=1$ computed using a $\mathbb P_4$ FR approximation on a uniform $800 \times 400$ mesh are shown in \cref{fig:shocktube_re} for varying values of the dynamic viscosity $\mu$. As the Reynolds number was increased, progressively smaller-scale features became evident in the flow, particularly along the contact line and in the separation region. The discontinuities as well as the small vortical structures in the flow were both well-resolved. The predicted results show relatively good agreement with the results of \citet{Guermond2021}, although some small-scale features at higher Reynolds numbers were notably more pronounced in the present work. These results suggest that the entropy filter used in conjunction with an operator splitting approach is an effective tool for approximating mixed hyperbolic-parabolic systems that exhibit discontinuities as well as complex small-scale flow structures. 

    \begin{figure}[htbp!]
        \centering 
        ~~~~~~
        \adjustbox{width=0.7\linewidth,valign=b}{\includegraphics[width=\textwidth]{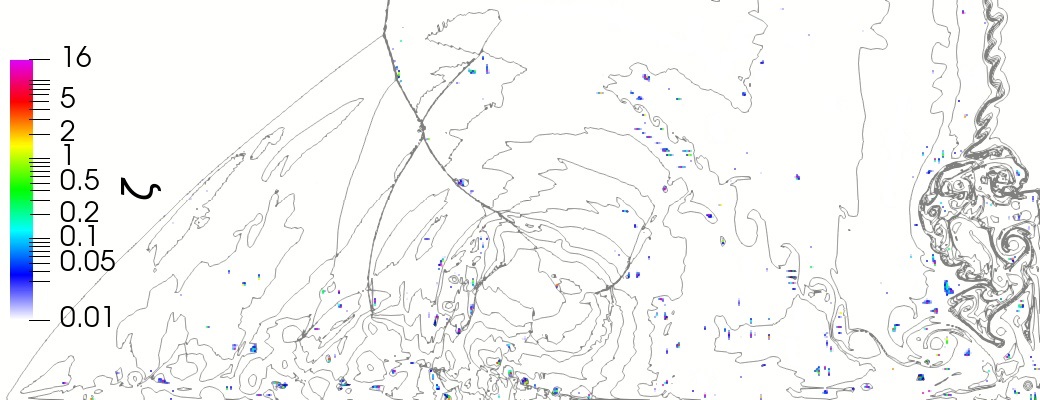}}
        \newline
        \caption{\label{fig:shocktube_zeta} Distribution of the filter parameter $\zeta$ overlaid on isocontours of density for the case of \cref{fig:st_10k}.}
    \end{figure}
    
To investigate the behavior of the filter for viscous flows, the distribution of the filter parameter $\zeta$ overlaid on isocontours of density is shown in \cref{fig:shocktube_zeta} for the case presented in \cref{fig:st_10k}. The distribution for the present case was more sporadic than in the inviscid cases, with more frequent activation near shock-boundary layer interactions (some of which are more easily observed with the Schlieren diagram in \cref{fig:st_10k} than with the isocontours in \cref{fig:shocktube_zeta}). It was also observed that for a small portion of the elements within the mesh, the addition of the viscous component caused the initially positivity-preserving solution to violate the positivity constraints. Approximately 16 elements required additional filtering, but the amount of filtering required was very minor, with values of $\zeta$ on the order of $10^{-4}$.

\subsubsection{Transonic Delta Wing}
As a final evaluation of the proposed approach for complex aeronautical applications including three-dimensional high Reynolds number flows computed on unstructured meshes, the test case of large eddy simulation around a transonic VFE-2 delta wing was considered. The geometry, introduced by \citet{Chu1996}, consists of a sharp leading edge delta wing with a sweep angle of $65^{\circ}$ and a thickness to root chord ratio of $3.4\%$. To match the experimental setup of \citet{Konrath2006}, the Reynolds number is set to $3{\cdot}10^6$ based on a unit root chord (i.e., $2{\cdot}10^6$ based on the mean aerodynamic chord), and the freestream Mach number is set to 0.8, yielding locally supersonic turbulent flow on the suction side of the wing. For these geometry and flow conditions, the separation and roll-up of the primary vortex is triggered by the sharp leading edge which immediately transitions into turbulence. Secondary and tertiary vortices are then subsequently formed through near-wall viscous interactions. An angle of attack of $20.5^{\circ}$ was chosen as these conditions yield strongly nonlinear behavior on the suction side as the primary vortex is on the precipice of vortex burst.

    \begin{figure}[htbp!]
        \centering
        \adjustbox{width=0.5\linewidth,valign=b}{\includegraphics[width=\textwidth]{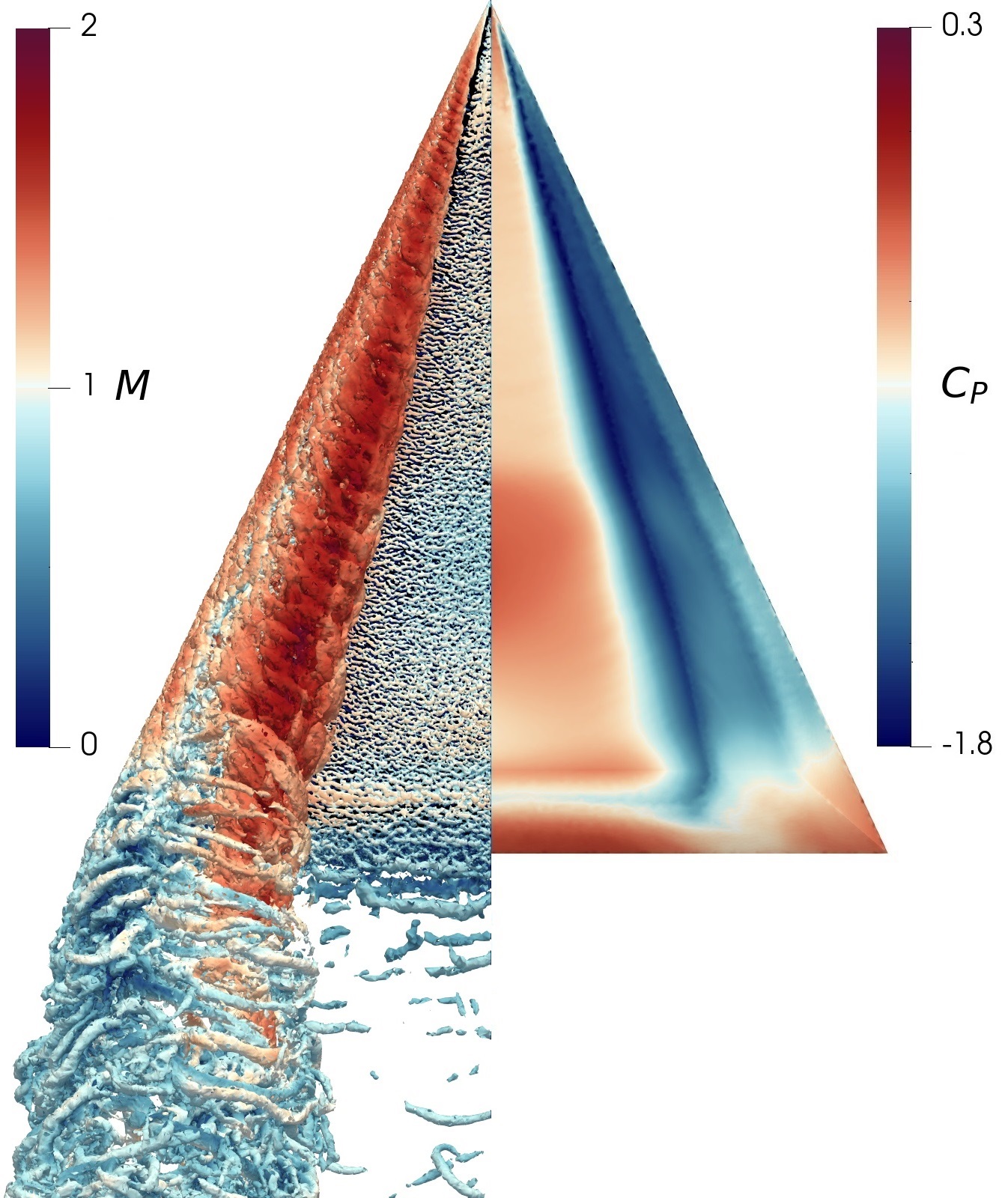}}
        \caption{\label{fig:delta} Isosurface of instantaneous Q-criterion colored by local Mach number (left) and time-averaged surface pressure coefficient contours (right) for the transonic VFE-2 delta wing computed using a $\mathbb P_3$ FR approximation.}
    \end{figure}
    
A second-order unstructured mesh was generated by extruding a triangular surface mesh in the near wall region and filling the remaining domain with tetrahedral elements. The resolution of the surface mesh was approximately uniform with an average edge length of $h = 0.003$, and the near-wall prismatic region was extruded over a length of $0.02$. To simplify the meshing, the sting and fairing from the experimental model were neglected. Since the separation and roll-up of the primary vortex is primarily governed by inviscid effects due to the sharp leading edge \citep{Luckring2002}, the wall-normal resolution requirements are relaxed in comparison to standard wall-resolved large eddy simulation \citep{Dzanic2019}. As such, the wall-normal resolution was set such that the $y^+$ value of the first solution point away from the wall was approximately 30 (i.e., the element $y^+$ was ${\sim}120$) based on equivalent flat plate conditions matching the freestream flow and root chord. In the separated region above the suction side, the resolution of the tetrahedral mesh was such that the average edge length was approximately $h = 0.01$. The overall mesh consisted of 4.1 million elements, yielding a total of 104 million degrees of freedom per component with a $\mathbb P_3$ approximation. At the wing, no slip adiabatic wall boundary conditions were applied, and characteristic Riemann invariant boundary conditions were used for the farfield.  

At startup, the approximation order and Reynolds number were progressively increased until the initial transients were convected away. Averaging was then performed for a time period corresponding to five flows over root chord. The results of a $\mathbb P_3$ FR approximation are presented in \cref{fig:delta}, showing instantaneous Q-criterion isosurfaces colored by local Mach number and time-averaged surface pressure coefficient contours. The Q-criterion isosurfaces showed the canonical vortex structure of delta wing flows with the primary and secondary vortices evident in the flow. Due to the relatively high angle of attack, a local Mach number of approximately 2 was observed near the primary vortex core, and indication of the onset of vortex burst was seen near the trailing edge. Furthermore, the average surface pressure coefficient contours clearly showed the presence of a primary and secondary vortex, and good agreement was observed between the predicted results and the experimental results of \citet{Konrath2006}.

\begin{figure}[htbp!]
        \centering
        \subfloat[$x/c = 0.2$]{\adjustbox{width=0.33\linewidth,valign=b}{\input{figs/delta_x2}}}
        ~
        \subfloat[$x/c = 0.4$]{\adjustbox{width=0.33\linewidth,valign=b}{\input{figs/delta_x4}}}
        \newline
        \subfloat[$x/c = 0.6$]{\adjustbox{width=0.33\linewidth,valign=b}{\input{figs/delta_x6}}}
        ~
        \subfloat[$x/c = 0.8$]{\adjustbox{width=0.33\linewidth,valign=b}{\input{figs/delta_x8}}}
        \newline
        \caption{\label{fig:delta_cp} Surface pressure coefficient at various streamwise locations for the transonic VFE-2 delta wing computed using a $\mathbb P_3$ FR approximation. Experimental results of \citet{Konrath2006} shown for reference. Dotted line denotes the sonic pressure coefficient. Spanwise extent is normalized by the local semispan $s$. }
    \end{figure}
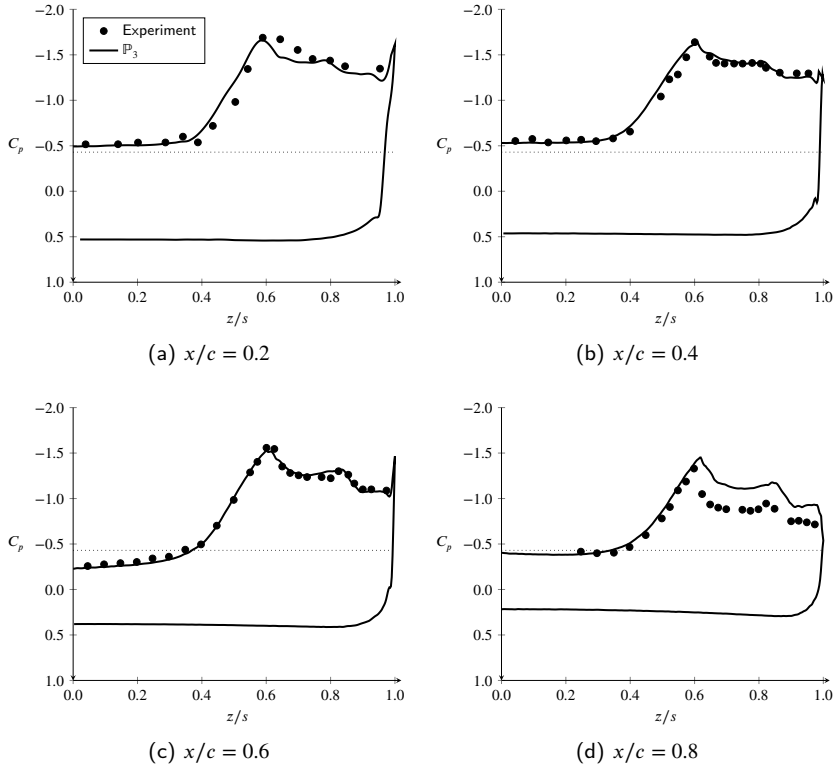

The surface pressure coefficient at various streamwise locations is presented in \cref{fig:delta_cp} in comparison to the experimental data. At $x/c = 0.2$, $0.4$, and $0.6$, excellent agreement was observed, both in terms of the location and the magnitude of the primary and secondary vortex pressure peaks. At $x/c = 0.8$, the location and magnitude of the primary vortex pressure peak was well-resolved, but a slight overprediction in the magnitude of the secondary vortex pressure peak was seen. This effect may be attributed to the strong sensitivity of the flow regime prior to vortex burst in conjunction with the slight discrepancies between the experimental and computational geometry at the trailing edge of the wing. Overall, the computational results were in-line with the experimental data, indicating that the entropy filtering approach can be successfully used for complex supersonic turbulent flows. 
 
%% file: figs/sod_p3.tex
    \begin{tikzpicture}[spy using outlines={rectangle, height=3cm,width=2.3cm, magnification=3, connect spies}]
		\begin{axis}[name=plot1,
		    axis line style={latex-latex},
		    axis x line=left,
            axis y line=left,
            clip mode=individual,
		    xlabel={$x$},
		    xtick={0,0.2,0.4,0.6,0.8,1},
    		xmin=0,
    		xmax=1,
    		x tick label style={
        		/pgf/number format/.cd,
            	fixed,
            	fixed zerofill,
            	precision=1,
        	    /tikz/.cd},
    		ylabel={$\rho$},
    		ylabel style={rotate=-90},
    		ytick={0,0.2,0.4,0.6,0.8,1},
    		ymin=0,
    		ymax=1.05,
    		y tick label style={
        		/pgf/number format/.cd,
            	fixed,
            	fixed zerofill,
            	precision=1,
        	    /tikz/.cd},
    		legend style={at={(0.03,0.03)},anchor=south west,font=\small},
    		legend cell align={left},
    		style={font=\normalsize}]
    		
    		\addplot[color=gray, style={ultra thin}, only marks, mark=o, mark options={scale=0.7}, mark repeat = 5, mark phase =0]
				table[x=x,y=r,col sep=comma,unbounded coords=jump]{./figs/data/sod_exact_reference.csv};
    		\addlegendentry{Exact}
    		
    		\addplot[color=red!80!black, style={very thick}]
				table[x=x,y=r,col sep=comma,unbounded coords=jump]{./figs/data/sod_p3_200dof.csv};
			\addlegendentry{$\mathbb P_3$}
			
		\end{axis}

	\end{tikzpicture}

%% file: figs/sod_p5.tex
    \begin{tikzpicture}[spy using outlines={rectangle, height=3cm,width=2.3cm, magnification=3, connect spies}]
		\begin{axis}[name=plot1,
		    axis line style={latex-latex},
		    axis x line=left,
            axis y line=left,
            clip mode=individual,
		    xlabel={$x$},
		    xtick={0,0.2,0.4,0.6,0.8,1},
    		xmin=0,
    		xmax=1,
    		x tick label style={
        		/pgf/number format/.cd,
            	fixed,
            	fixed zerofill,
            	precision=1,
        	    /tikz/.cd},
    		ylabel={$\rho$},
    		ylabel style={rotate=-90},
    		ytick={0,0.2,0.4,0.6,0.8,1},
    		ymin=0,
    		ymax=1.05,
    		y tick label style={
        		/pgf/number format/.cd,
            	fixed,
            	fixed zerofill,
            	precision=1,
        	    /tikz/.cd},
    		legend style={at={(0.03,0.03)},anchor=south west,font=\small},
    		legend cell align={left},
    		style={font=\normalsize}]
    		
    		\addplot[color=gray, style={ultra thin}, only marks, mark=o, mark options={scale=0.7}, mark repeat = 5, mark phase =0]
				table[x=x,y=r,col sep=comma,unbounded coords=jump]{./figs/data/sod_exact_reference.csv};
    		\addlegendentry{Exact}
    		
    		\addplot[color=red!80!black, style={very thick}]
				table[x=x,y=r,col sep=comma,unbounded coords=jump]{./figs/data/sod_p5_200dof.csv};
			\addlegendentry{$\mathbb P_5$}
			
		\end{axis}

	\end{tikzpicture}

%% file: figs/shuosher100.tex
     \begin{tikzpicture}[spy using outlines={rectangle, height=3cm,width=2.3cm, magnification=3, connect spies}]
		\begin{axis}[name=plot1,
		    axis line style={latex-latex},
		    axis x line=left,
            axis y line=left,
		    clip mode=individual,
		    xlabel={$x$},
		    xtick={-5,-2.5,0,2.5,5},
    		xmin=-5,
    		xmax=5,
    		x tick label style={
        		/pgf/number format/.cd,
            	fixed,
            	fixed zerofill,
            	precision=1,
        	    /tikz/.cd},
    		ylabel={$\rho$},
    		ylabel style={rotate=-90},
    		ytick={0,1,2,3,4,5},
    		ymin=0,
    		ymax=5,
    		y tick label style={
        		/pgf/number format/.cd,
            	fixed,
            	fixed zerofill,
            	precision=0,
        	    /tikz/.cd},
    		legend style={at={(0.03,0.075)},anchor=south west,font=\small},
    		legend cell align={left},
    		style={font=\normalsize}]
    		
    		\addplot[color=gray, style={ultra thin}, only marks, mark=o, mark options={scale=0.5}, mark repeat = 2, mark phase =0]
				table[x=x,y=d,col sep=comma,unbounded coords=jump]{./figs/data/osher_p0_2000.csv};
    		\addlegendentry{Reference}
    		
    		\addplot[color=red!80!black, style={very thick}]
				table[x expr={\thisrow{x}-5},y=r,col sep=comma,unbounded coords=jump]{./figs/data/shuosher_p3_n100.csv};
    		\addlegendentry{$\mathbb P_3$}
    		
% 			\addplot[color=black, style={ultra thin}, only marks, mark=o, mark options={scale=0.4}, mark repeat = 3, mark phase =0]
% 				table[x expr={\thisrow{x}-5},y=r,col sep=comma,unbounded coords=jump]{./figs/data/shuosher_p3_n200.csv};
    % 		\addlegendentry{$N = 200$}

		\end{axis}

	\end{tikzpicture}

%% file: figs/shuosher200.tex
     \begin{tikzpicture}[spy using outlines={rectangle, height=3cm,width=2.3cm, magnification=3, connect spies}]
		\begin{axis}[name=plot1,
		    axis line style={latex-latex},
		    axis x line=left,
            axis y line=left,
		    clip mode=individual,
		    xlabel={$x$},
		    xtick={-5,-2.5,0,2.5,5},
    		xmin=-5,
    		xmax=5,
    		x tick label style={
        		/pgf/number format/.cd,
            	fixed,
            	fixed zerofill,
            	precision=1,
        	    /tikz/.cd},
    		ylabel={$\rho$},
    		ylabel style={rotate=-90},
    		ytick={0,1,2,3,4,5},
    		ymin=0,
    		ymax=5,
    		y tick label style={
        		/pgf/number format/.cd,
            	fixed,
            	fixed zerofill,
            	precision=0,
        	    /tikz/.cd},
    		legend style={at={(0.03,0.075)},anchor=south west,font=\small},
    		legend cell align={left},
    		style={font=\normalsize}]
    		
    		\addplot[color=gray, style={ultra thin}, only marks, mark=o, mark options={scale=0.5}, mark repeat = 2, mark phase =0]
				table[x=x,y=d,col sep=comma,unbounded coords=jump]{./figs/data/osher_p0_2000.csv};
    % 		\addlegendentry{Reference}

% 			\addplot[color=black, style={ultra thin}, only marks, mark=+, mark options={scale=0.4}, mark repeat = 3, mark phase =0]
% 				table[x expr={\thisrow{x}-5},y=r,col sep=comma,unbounded coords=jump]{./figs/data/shuosher_p3_n100.csv};
    % 		\addlegendentry{$\mathbb P_3$}
    		
    		\addplot[color=red!80!black, style={thick}]
				table[x expr={\thisrow{x}-5},y=r,col sep=comma,unbounded coords=jump,each nth point=2]{./figs/data/shuosher_p3_n200.csv};
    % 		\addlegendentry{$N = 200$}

		\end{axis}

	\end{tikzpicture}

%% file: figs/tgv_dedt.tex
\begin{tikzpicture}[spy using outlines={rectangle, height=3cm,width=2.5cm, magnification=3, connect spies}]
    \begin{axis}
    [
        axis line style={latex-latex},
        axis y line=left,
        axis x line=left,
        clip mode=individual,
        xmode=linear, % not log
        ymode=linear, % not log
        xlabel = {$t$},
        ylabel = {$\varepsilon_K$, \textcolor{red!80!black}{$\varepsilon_E$}},
        xmin = 0, xmax = 20,
        ymin = 0.00, ymax = 0.014,
        legend cell align={left},
        legend style={font=\scriptsize, at={(1.0, 1.0)}, anchor=north east},
        ytick = {0,0.004,0.008,0.012,0.016},
        x tick label style={/pgf/number format/.cd, fixed, fixed zerofill, precision=0, /tikz/.cd},
        y tick label style={/pgf/number format/.cd, fixed, fixed zerofill, precision=1, /tikz/.cd},	
        scale = 0.9
    ]
        
        \addplot[ color=gray, style={thick}, only marks, mark=o, mark options={scale=0.5}, mark repeat = 20, mark phase = 0] table[x=t, y=dkep, col sep=comma, mark=*]{./figs/data/van_rees_2011_dkep_PSP512.csv};
        \addlegendentry{DNS};
        
        \addplot[color=black, style={dotted, thick}] table[x=t, y=dEdt, col sep=comma]{./figs/data/tgvstats96dofp3filtere4_post.csv};
        \addlegendentry{$96^3$ DoF};
        
        \addplot[color=black, style={dashed, thick}] table[x=t, y=dEdt, col sep=comma]{./figs/data/tgvstats128dofp3filtere4_post.csv};
        \addlegendentry{$128^3$ DoF};
        
        \addplot[color=black, style={thick}] table[x=t, y=dEdt, col sep=comma]{./figs/data/tgvstats160dofp3filtere4v2_post.csv};
        \addlegendentry{$160^3$ DoF };

        \addplot[color=red!80!black, style={dotted, thick}] table[x=t, y=enst, col sep=comma]{./figs/data/tgvstats96dofp3filtere4_post.csv};
        % \addlegendentry{$96^3$ DoF};
        
        \addplot[color=red!80!black, style={dashed, thick}] table[x=t, y=enst, col sep=comma]{./figs/data/tgvstats128dofp3filtere4_post.csv};
        % \addlegendentry{$128^3$ DoF};
        
        \addplot[color=red!80!black, style={thick}] table[x=t, y=enst, col sep=comma]{./figs/data/tgvstats160dofp3filtere4v2_post.csv};
        % \addlegendentry{$160^3$ DoF };
        
    \end{axis}

\end{tikzpicture}

%% file: figs/tgv_enst.tex
\begin{tikzpicture}[spy using outlines={rectangle, height=3cm,width=2.5cm, magnification=3, connect spies}]
    \begin{axis}
    [
        axis line style={latex-latex},
        axis y line=left,
        axis x line=left,
        clip mode=individual,
        xmode=linear, % not log
        ymode=linear, % not log
        xlabel = {$t$},
        ylabel = {$\varepsilon_K$, \textcolor{red!80!black}{$\varepsilon_E$}},
        xmin = 0, xmax = 20,
        ymin = 0.00, ymax = 0.014,
        legend cell align={left},
        legend style={font=\scriptsize, at={(1.0, 1.0)}, anchor=north east},
        ytick = {0,0.004,0.008,0.012,0.016},
        x tick label style={/pgf/number format/.cd, fixed, fixed zerofill, precision=0, /tikz/.cd},
        y tick label style={/pgf/number format/.cd, fixed, fixed zerofill, precision=1, /tikz/.cd},	
        scale = 0.9
    ]
        
        \addplot[ color=gray, style={thick}, only marks, mark=o, mark options={scale=0.5}, mark repeat = 20, mark phase = 0] table[x=t, y=dkep, col sep=comma, mark=*]{./figs/data/van_rees_2011_dkep_PSP512.csv};
        \addlegendentry{DNS};
        
        \addplot[color=black, style={dotted, thick}] table[x=t, y=dEdt, col sep=comma]{./figs/data/tgvstats128dofp3filtere4_post.csv};
        \addlegendentry{$\mathbb P_3$};
        
        \addplot[color=black, style={dashed, thick}] table[x=t, y=dEdt, col sep=comma]{./figs/data/tgvstats120dofp4filtere4v2_post.csv};
        \addlegendentry{$\mathbb P_4$};
        
        \addplot[color=black, style={thick}] table[x=t, y=dEdt, col sep=comma]{./figs/data/tgvstats120dofp5filtere4v2_post.csv};
        \addlegendentry{$\mathbb P_5$};
        
        \addplot[color=red!80!black, style={dotted, thick}] table[x=t, y=enst, col sep=comma]{./figs/data/tgvstats128dofp3filtere4_post.csv};
        % \addlegendentry{$96^3$ DoF};
        
        \addplot[color=red!80!black, style={dashed, thick}] table[x=t, y=enst, col sep=comma]{./figs/data/tgvstats120dofp4filtere4v2_post.csv};
        % \addlegendentry{$128^3$ DoF};
        
        \addplot[color=red!80!black, style={thick}] table[x=t, y=enst, col sep=comma]{./figs/data/tgvstats120dofp5filtere4v2_post.csv};
        
    \end{axis}

\end{tikzpicture}

%% file: figs/delta_x2.tex
    \begin{tikzpicture}[spy using outlines={rectangle, height=3cm,width=2.3cm, magnification=3, connect spies}]
		\begin{axis}[name=plot1,
		    axis line style={latex-latex},
		    axis x line=left,
            axis y line=left,
            clip mode=individual,
		    xlabel={$z/s$},
		    xtick={0,0.2,0.4,0.6,0.8,1},
    		xmin=0,
    		xmax=1.02,
    	    x tick label style={
        		/pgf/number format/.cd,
            	fixed,
            	fixed zerofill,
            	precision=1,
        	    /tikz/.cd},
    		ylabel={$C_p$},
    		ylabel style={rotate=-90},
    		ytick={1.0, 0.5, 0.0, -0.5, -1.0, -1.5, -2.0},
    		ymin=-2.0,
    		ymax=1,
    		y dir=reverse,
    		y tick label style={
        		/pgf/number format/.cd,
            	fixed,
            	fixed zerofill,
            	precision=1,
        	    /tikz/.cd},
    		legend style={at={(0.03, 0.97)},anchor=north west},
    		legend cell align={left}]
    		
			\addplot[color=black, style={thin}, only marks, mark=*, mark options={scale=1.0}]
				table[x=z,y=cp,col sep=comma,unbounded coords=jump]{./figs/data/vfe2_ref_x0p2.csv};
    		\addlegendentry{Experiment}
    		
			\addplot[color=black, style={very thick}]
				table[x=z,y=cp,col sep=comma,unbounded coords=jump]{./figs/data/delta_x0p2.csv};
    		\addlegendentry{$\mathbb P_3$}
    		
        \draw [dotted] (0,-0.43) -- (1,-0.43);

		\end{axis} 		
	\end{tikzpicture}

%% file: figs/delta_x4.tex
    \begin{tikzpicture}[spy using outlines={rectangle, height=3cm,width=2.3cm, magnification=3, connect spies}]
		\begin{axis}[name=plot1,
		    axis line style={latex-latex},
		    axis x line=left,
            axis y line=left,
            clip mode=individual,
		    xlabel={$z/s$},
		    xtick={0,0.2,0.4,0.6,0.8,1},
    		xmin=0,
    		xmax=1.02,
    	    x tick label style={
        		/pgf/number format/.cd,
            	fixed,
            	fixed zerofill,
            	precision=1,
        	    /tikz/.cd},
    		ylabel={$C_p$},
    		ylabel style={rotate=-90},
    		ytick={1.0, 0.5, 0.0, -0.5, -1.0, -1.5, -2.0},
    		ymin=-2.0,
    		ymax=1,
    		y dir=reverse,
    		y tick label style={
        		/pgf/number format/.cd,
            	fixed,
            	fixed zerofill,
            	precision=1,
        	    /tikz/.cd},
    		legend style={at={(0.03, 0.97)},anchor=north west},
    		legend cell align={left}]
    		
			\addplot[color=black, style={thin}, only marks, mark=*, mark options={scale=1.0}]
				table[x=z,y=cp,col sep=comma,unbounded coords=jump]{./figs/data/vfe2_ref_x0p4.csv};
    % 		\addlegendentry{Experiment}
    		
			\addplot[color=black, style={very thick}]
				table[x=z,y=cp,col sep=comma,unbounded coords=jump]{./figs/data/delta_x0p4.csv};
    % 		\addlegendentry{$\mathbb P_3$}
    		
        \draw [dotted] (0,-0.43) -- (1,-0.43);

		\end{axis} 		
	\end{tikzpicture}

%% file: figs/delta_x6.tex
    \begin{tikzpicture}[spy using outlines={rectangle, height=3cm,width=2.3cm, magnification=3, connect spies}]
		\begin{axis}[name=plot1,
		    axis line style={latex-latex},
		    axis x line=left,
            axis y line=left,
            clip mode=individual,
		    xlabel={$z/s$},
		    xtick={0,0.2,0.4,0.6,0.8,1},
    		xmin=0,
    		xmax=1.02,
    	    x tick label style={
        		/pgf/number format/.cd,
            	fixed,
            	fixed zerofill,
            	precision=1,
        	    /tikz/.cd},
    		ylabel={$C_p$},
    		ylabel style={rotate=-90},
    		ytick={1.0, 0.5, 0.0, -0.5, -1.0, -1.5, -2.0},
    		ymin=-2.0,
    		ymax=1,
    		y dir=reverse,
    		y tick label style={
        		/pgf/number format/.cd,
            	fixed,
            	fixed zerofill,
            	precision=1,
        	    /tikz/.cd},
    		legend style={at={(0.03, 0.97)},anchor=north west},
    		legend cell align={left}]
    		
			\addplot[color=black, style={thin}, only marks, mark=*, mark options={scale=1.0}]
				table[x=z,y=cp,col sep=comma,unbounded coords=jump]{./figs/data/vfe2_ref_x0p6.csv};
    % 		\addlegendentry{Experiment}
    		
			\addplot[color=black, style={very thick}]
				table[x=z,y=cp,col sep=comma,unbounded coords=jump]{./figs/data/delta_x0p6.csv};
    % 		\addlegendentry{$\mathbb P_3$}
    		
        \draw [dotted] (0,-0.43) -- (1,-0.43);

		\end{axis} 		
	\end{tikzpicture}

%% file: figs/delta_x8.tex
    \begin{tikzpicture}[spy using outlines={rectangle, height=3cm,width=2.3cm, magnification=3, connect spies}]
		\begin{axis}[name=plot1,
		    axis line style={latex-latex},
		    axis x line=left,
            axis y line=left,
            clip mode=individual,
		    xlabel={$z/s$},
		    xtick={0,0.2,0.4,0.6,0.8,1},
    		xmin=0,
    		xmax=1.02,
    	    x tick label style={
        		/pgf/number format/.cd,
            	fixed,
            	fixed zerofill,
            	precision=1,
        	    /tikz/.cd},
    		ylabel={$C_p$},
    		ylabel style={rotate=-90},
    		ytick={1.0, 0.5, 0.0, -0.5, -1.0, -1.5, -2.0},
    		ymin=-2.0,
    		ymax=1,
    		y dir=reverse,
    		y tick label style={
        		/pgf/number format/.cd,
            	fixed,
            	fixed zerofill,
            	precision=1,
        	    /tikz/.cd},
    		legend style={at={(0.03, 0.97)},anchor=north west},
    		legend cell align={left}]
    		
			\addplot[color=black, style={thin}, only marks, mark=*, mark options={scale=1.0}]
				table[x=z,y=cp,col sep=comma,unbounded coords=jump]{./figs/data/vfe2_ref_x0p8.csv};
    % 		\addlegendentry{Experiment}
    		
			\addplot[color=black, style={very thick}]
				table[x=z,y=cp,col sep=comma,unbounded coords=jump]{./figs/data/delta_x0p8.csv};
    % 		\addlegendentry{$\mathbb P_3$}
    		
        \draw [dotted] (0,-0.43) -- (1,-0.43);
    		
		\end{axis} 		
	\end{tikzpicture}

%% file: conclusions.tex
\section{Conclusions}\label{sec:conclusion}
In this work, a novel adaptive filtering approach for shock capturing in discontinuous spectral element methods is presented. The proposed entropy filtering approach formulates the filtering operation as an element-wise scalar optimization problem designed to enforce constraints such as positivity and a local discrete minimum entropy principle. As a result, this method is free of problem-dependent tunable parameters and can be easily implemented on arbitrary unstructured meshes with minimal computational effort. To evaluate the efficacy of the proposed filtering approach, a series of numerical experiments was performed on multi-dimensional hyperbolic and mixed hyperbolic-parabolic conservation laws such as the Euler and Navier--Stokes equations. Even for extremely strong shocks, sub-element resolution of the discontinuities was generally obtained with minimal spurious oscillations, and high-order accuracy was recovered for smooth solutions. Furthermore, for more complex problems including shock-vortex interactions and compressible turbulent flows, the filter was able to robustly resolve discontinuities without unnecessarily dissipating small-scale vortical structures in the flow. These results indicate that the proposed entropy filtering approach can be an effective and robust technique for shock capturing in discontinuous spectral element approximations of hyperbolic and mixed hyperbolic-parabolic conservation laws. In future work, the proposed approach can be improved by optimizing the calculation of the filter strength through more sophisticated root bracketing methods and stopping criteria, increasing the locality of the filtering kernel to mimic that of methods such as artificial viscosity and subcell approaches, modifying the approach to remove the need for operator splitting in mixed systems, and applying a more robust method for calculating the entropy tolerance.

%% file: app1.tex
\section{Algorithmic Details}\label{app:algo}
An in-depth description of the computational details of the proposed filtering approach is presented below. In \cref{alg:entropy}, the computation of the entropy constraints is presented. This method computes the discrete local minima of the entropy function $\sigma_*^k$ at the solution nodes for each element $\Omega_k$. The entropy constraint $\sigma_{\min}^k$ is then taken as the minima of $\sigma_*^k$ across the face-neighbors of $\Omega_k$, including $\Omega_k$ itself.

\begin{algorithm}
\caption{Calculate entropy constraints}
{INPUT}: $\mathbf{u}$ // Global solution \newline
{OUTPUT}: $\sigma_{\min}$ // Global entropy constraints\newline

GetEntropyConstraints($\mathbf{u}$) :
\begin{algorithmic} \label{alg:entropy}
\FOR{$\Omega_k \in \Omega$}
\item[]$\sigma_{*}^k = \underset{i \in S(\Omega_k)}{\min} \sigma(\mathbf{u}_i)$
\ENDFOR
\FOR{$\Omega_k \in \Omega$}
\item[]$\sigma_{\min}^k = \underset{j \in \mathcal A (k)}{\min} \sigma_{*}^j$
\ENDFOR
\item[] return $\sigma_{\min}$
\end{algorithmic}
\end{algorithm}

In \cref{alg:filter}, the element-wise filtering operation is outlined. This method takes in the unfiltered element-wise solution $\mathbf{u}$ and its associated entropy constraint $\sigma_{\min}$ as the input in addition to a precomputed Vandermonde matrix $\mathbf{V}$ (from the nodal basis to a modal basis). If the discrete solution is within bounds, it returns the unfiltered solution. If not, it performs $n_{\text{iters}}$ iterations of the bisection approach. The initial guess for the upper bound of the bisection method is taken as $-\log(\varepsilon)$ for some value of machine precision $\varepsilon$ as this effectively approximates an infinite filter strength down to machine precision.
\begin{algorithm}
\caption{Filter solution}
{INPUT}: $\mathbf{u}$, $\sigma_{\min}$ // Element-wise solution, element-wise entropy constraint\newline
{OUTPUT}: $\widetilde{\mathbf{u}}$ // Element-wise filtered solution\newline

FilterSolution($\mathbf{u}$, $\sigma_{\min}$) :
\begin{algorithmic}\label{alg:filter}

\item[] // PARAMETERS
\item[] $\rho_{\min} = 10^{-8}$ // Minimum density
\item[] $P_{\min} = 10^{-8}$ // Minimum pressure
\item[]$\epsilon_\sigma = 10^{-4}$ // Entropy tolerance
\item[]$\varepsilon = 10^{-8}$ // Machine precision  (FP32)
\item[]$n_{\text{iters}} = 20$ // Number of filter iterations
\item[]
\item[]// Compute minima within element
\item[]$[\rho_{*}, P_{*}, \sigma_{*}] = \underset{i \in S}{\min} [\rho(\mathbf{u}_i), P(\mathbf{u}_i), \sigma(\mathbf{u}_i)]$
\item[]
\item[]// Return unfiltered solution if within bounds 
\IF{($\rho_* \geq \rho_{\min}$) \& ($P_* \geq P_{\min}$) \& ($\sigma_* \geq \sigma_{\min} - \epsilon_\sigma$)}
  \STATE return $\mathbf{u}$
\ELSE
\item[]$\zeta_1$ = 0.0
\item[]$\zeta_2$ = $-\log(\varepsilon)$
\item[]
\item[]// Compute modal basis
\item[]$\hat{\mathbf{u}} = \mathbf{V}^{-1} \mathbf{u}$ 
\item[]
\FOR{$i\in\{0,..,n_{\text{iters}}\}$}
\item[]$\zeta_3 = 0.5(\zeta_1 + \zeta_2)$
\item[]$\widetilde{\mathbf{u}} = \mathbf{V} \left(\hat{\mathbf{u}} \odot \exp (-\zeta_3 \mathbf{p}^2) \right)$ 
\item[]
\item[]// Compute minima within element
\item[]$[\rho_{*}, P_{*}, \sigma_{*}] = \underset{i \in S}{\min} [\rho(\widetilde{\mathbf{u}}_i), P(\widetilde{\mathbf{u}}_i), \sigma(\widetilde{\mathbf{u}}_i)]$
\item[]
\item[] // Choose new midpoint
\IF{($\rho_* \geq \rho_{\min}$) \& ($P_* \geq P_{\min}$) \& ($\sigma_* \geq \sigma_{\min} - \epsilon_\sigma$)}
\item[]$\zeta_2 = \zeta_3$ 
\ELSE
\item[]$\zeta_1 = \zeta_3$ 
\ENDIF
\ENDFOR
\item[]
\item[] // Compute filtered solution with bounds-preserving filter 
\item[]$\zeta = \zeta_2$
\item[]$\widetilde{\mathbf{u}} = \mathbf{V} \left(\hat{\mathbf{u}} \odot \exp (-\zeta \mathbf{p}^2) \right)$ 
\item[]return $\widetilde{\mathbf{u}}$
\item[]
\ENDIF
\item[]
\end{algorithmic}
\end{algorithm}

In \cref{alg:stepper_hyp}, an outline of the temporal integration approach utilizing an SSP-RK3 scheme \emph {for purely hyperbolic systems} is given. For a given time step $\Delta t$ and number of time steps $nt$, the method computes each intermediate stage of the temporal integration scheme and filters the solution based on the constraints computed on the previous stage. A similar outline is shown in \cref{alg:stepper_mix} for mixed hyperbolic-parabolic systems that details the operator splitting approach. 

\begin{algorithm}
\caption{SSP-RK3 stepper for hyperbolic systems}
// PARAMETERS\newline
$\Delta t$ = ... // Time step\newline
$n_t$ = ... // Number of time steps \newline
\newline
// Step in time 
\begin{algorithmic} \label{alg:stepper_hyp}
\FOR{$n\in \{0,..,n_t\}$} 
\item[]// Calculate entropy constraints for each element
\item[]$\sigma_{\min} =$ GetEntropyConstraints($\mathbf{u}^n$) 
\item[]// Compute first stage
\item[]$ \mathbf{u}^{(1)} = \mathbf{u}^n + \Delta t \left(-\boldsymbol{\nabla} {\cdot}\mathbf{F} (\mathbf{u}^n) \right )$ 
\item[]// Check if constraints are satisfied and filter if not
\FOR{$\Omega_k \in \Omega$} 
\item[]$\mathbf{u}^{(1)}_k$ = FilterSolution($\mathbf{u}^{(1)}_k$, $\sigma_{\min}^k$)
\ENDFOR
\item[]
\item[] // Repeat procedure for second stage
\item[]$\sigma_{\min} =$ GetEntropyConstraints($\mathbf{u}^{(1)}$) 
\item[]$ \mathbf{u}^{(2)} = \frac{3}{4}\mathbf{u}^n + \frac{1}{4}\mathbf{u}^{(1)} + \Delta t \left(-\boldsymbol{\nabla} {\cdot} \mathbf{F} (\mathbf{u}^{(1)}) \right )$ 
\FOR{$\Omega_k \in \Omega$} 
\item[]$\mathbf{u}^{(2)}_k$ = FilterSolution($\mathbf{u}^{(2)}_k$, $\sigma_{\min}^k$)
\ENDFOR
\item[]
\item[] // Repeat procedure for third stage and compute next temporal step
\item[]$\sigma_{\min} =$ GetEntropyConstraints($\mathbf{u}^{(2)}$) 
\item[]$ \mathbf{u}^{n+1} = \frac{1}{3}\mathbf{u}^n + \frac{2}{3}\mathbf{u}^{(2)} + \Delta t \left(-\boldsymbol{\nabla} {\cdot} \mathbf{F} (\mathbf{u}^{(2)}) \right )$ 
\FOR{$\Omega_k \in \Omega$} 
\item[]$\mathbf{u}^{n+1}_k$ = FilterSolution($\mathbf{u}^{n+1}_k$, $\sigma_{\min}^k$)
\ENDFOR

\ENDFOR
\end{algorithmic}
\end{algorithm}

\begin{algorithm}
\caption{SSP-RK3 stepper for mixed hyperbolic-parabolic systems}
// PARAMETERS\newline
$\Delta t$ = ... // Time step\newline
$n_t$ = ... // Number of time steps \newline
\newline
// Step in time 
\begin{algorithmic} \label{alg:stepper_mix}
\FOR{$n\in\{0,..,nt\}$} 
\item[]// Calculate entropy constraints for each element
\item[]$\sigma_{\min} =$ GetEntropyConstraints($\mathbf{u}^n$) 
\item[]// Compute inviscid and viscous first stages
\item[]$ \mathbf{u}^{(1,I)} = \mathbf{u}^n + \Delta t \left(-\boldsymbol{\nabla} {\cdot}\mathbf{F}_I (\mathbf{u}^n) \right )$ 
\item[]$ \Delta \mathbf{u}^{(1,V)} = \Delta t \left(-\boldsymbol{\nabla} {\cdot}\mathbf{F}_V (\mathbf{u}^n) \right )$ 
\item[]// Check if constraints are satisfied and filter if not
\FOR{$\Omega_k \in \Omega$} 
\item[]// Enforce positivity/entropy constraints for inviscid component
\item[]$\mathbf{u}^{(1,I)}_k$ = FilterSolution($\mathbf{u}^{(1,I)}_k$, $\sigma_{\min}^k$)
\item[] // Add viscous component to filtered inviscid component
\item[] $\mathbf{u}^{(1)}_k = \mathbf{u}^{(1,I)}_k + \Delta\mathbf{u}^{(1,V)}_k$
\item[] // Enforce only positivity constraints for inviscid + viscous component
\item[]$\mathbf{u}^{(1)}_k$ = FilterSolution($\mathbf{u}^{(1)}_k$, $-\infty$)
\ENDFOR
\item[]
\item[] // Repeat procedure for second stage
\item[]$\sigma_{\min} =$ GetEntropyConstraints($\mathbf{u}^{(1)}$) 
\item[]$ \mathbf{u}^{(2, I)} = \frac{3}{4}\mathbf{u}^n + \frac{1}{4}\mathbf{u}^{(1)} + \Delta t \left(-\boldsymbol{\nabla} {\cdot} \mathbf{F}_I (\mathbf{u}^{(1)}) \right )$ 
\item[]$\Delta \mathbf{u}^{(2, V)} = \Delta t \left(-\boldsymbol{\nabla} {\cdot} \mathbf{F}_V (\mathbf{u}^{(1)}) \right )$ 
\FOR{$\Omega_k \in \Omega$} 
\item[]$\mathbf{u}^{(2,I)}_k$ = FilterSolution($\mathbf{u}^{(2,I)}_k$, $\sigma_{\min}^k$)
\item[] $\mathbf{u}^{(2)}_k = \mathbf{u}^{(2,I)}_k + \Delta\mathbf{u}^{(2,V)}_k$
\item[]$\mathbf{u}^{(2)}_k$ = FilterSolution($\mathbf{u}^{(2)}_k$, $-\infty$)
\ENDFOR
\item[]
\item[] // Repeat procedure for third stage and compute next temporal step
\item[]$\sigma_{\min} =$ GetEntropyConstraints($\mathbf{u}^{(2)}$) 
\item[]$ \mathbf{u}^{n+1, I} = \frac{1}{3}\mathbf{u}^n + \frac{2}{3}\mathbf{u}^{(2)} + \Delta t \left(-\boldsymbol{\nabla} {\cdot} \mathbf{F}_I (\mathbf{u}^{(2)}) \right )$ 
\item[]$\Delta \mathbf{u}^{n+1, V} = \Delta t \left(-\boldsymbol{\nabla} {\cdot} \mathbf{F}_V (\mathbf{u}^{(2)}) \right )$ 
\FOR{$\Omega_k \in \Omega$} 
\item[]$\mathbf{u}^{n+1, I}_k$ = FilterSolution($\mathbf{u}^{n+1, I}_k$, $\sigma_{\min}^k$)
\item[] $\mathbf{u}^{n+1}_k = \mathbf{u}^{n+1}_k + \Delta\mathbf{u}^{n+1, V}_k$
\item[]$\mathbf{u}^{n+1}_k$ = FilterSolution($\mathbf{u}^{n+1}_k$, $-\infty$)
\ENDFOR

\ENDFOR

\end{algorithmic}
\end{algorithm}